\documentclass[ejs,preprint]{imsart}

\RequirePackage[OT1]{fontenc}
\RequirePackage{amsthm,amsmath}
\usepackage[numbers,sort&compress]{natbib}
\RequirePackage[colorlinks,citecolor=blue,urlcolor=blue]{hyperref}
\RequirePackage{hypernat}

\RequirePackage{amssymb}
\RequirePackage{amsmath}
\RequirePackage{graphicx, rotating}
\RequirePackage{multirow}

\startlocaldefs
\numberwithin{equation}{section}
\theoremstyle{plain}

\newtheorem{theorem}{Theorem}[section]
\newtheorem{lemma}{Lemma}[section]
\newtheorem{proposition}{Proposition}[section]
\newtheorem{corollary}{Corollary}[section]

\newcommand{\bsb}{\boldsymbol}
\newcommand{\bsbX}{{\boldsymbol{X}}}
\newcommand{\bsbx}{{\boldsymbol{x}}}
\newcommand{\bsby}{{\boldsymbol{y}}}
\newcommand{\bsbb}{{\boldsymbol{\beta}}}
\newcommand{\bsbg}{{\boldsymbol{\gamma}}}

\newcommand{\bsbH}{{\boldsymbol{H}}}
\newcommand{\bsbI}{{\boldsymbol{I}}}

\newcommand{\bsbSig}{{\boldsymbol{\Sigma}}}

\newcommand{\bsbD}{{\boldsymbol{D}}}

\newcommand{\bsbU}{{\boldsymbol{U}}}

\newcommand{\bsba}{{\boldsymbol{\alpha}}}
\newcommand{\bsbA}{{\boldsymbol{A}}}

\newcommand{\bsbeps}{{\boldsymbol{\epsilon}}}

\newcommand{\bsbh}{{\boldsymbol{h}}}
\newcommand{\bsbS}{{\boldsymbol{S}}}

\newcommand{\gsgn}{{\widetilde{\mbox{sgn}}}}
\newcommand{\gSgn}{{\widetilde{\mbox{Sgn}}}}
\newcommand{\bsbs}{{\boldsymbol{s}}}

\endlocaldefs

\begin{document}

\begin{frontmatter}
\title{Thresholding-based Iterative Selection Procedures for Model Selection and Shrinkage}
\runtitle{TISP for Model Selection and Shrinkage}

\begin{aug}
\author{\fnms{Yiyuan} \snm{She}\ead[label=e1]{yshe@stat.fsu.edu}}
\runauthor{Y. She}

\affiliation{Florida State University}

\address{Department of Statistics\\
Florida State University\\
Tallahassee, FL 32306-4330\\
\printead{e1}\\
\phantom{E-mail:\ }}

\end{aug}

\begin{abstract}
This paper discusses a class of thresholding-based iterative selection procedures (TISP) for model selection and shrinkage. People have long before noticed the weakness of the convex $l_1$-constraint (or the soft-thresholding) in wavelets and have designed many different forms of nonconvex penalties to increase model sparsity and accuracy. But for a nonorthogonal regression matrix, there is great difficulty  in both investigating the  performance in theory and solving the problem  in computation. TISP provides a simple and efficient way to tackle this so that we successfully  borrow the rich results in the orthogonal design to solve the nonconvex penalized regression for a general design matrix.
Our starting point is, however, thresholding rules rather than penalty functions.
Indeed, there is  a universal connection between them.
But a drawback of the latter is its non-unique form, and our approach greatly facilitates the computation and the analysis. In fact, we are able to build the convergence theorem and explore  theoretical properties of the selection and  estimation via TISP nonasymptotically.
More importantly, a novel Hybrid-TISP is proposed based on hard-thresholding and ridge-thresholding. It provides a fusion between the $l_0$-penalty and the $l_2$-penalty, and adaptively achieves the right balance between shrinkage and selection in statistical modeling.
In practice, Hybrid-TISP shows superior performance in test-error and is parsimonious.

\end{abstract}



\end{frontmatter}

\section{Introduction}
Lasso~\cite{Tib} has attracted people's a lot of  attention recently because it provides an efficient and continuous way for variable selection, thereby achieving a stable sparse solution.
Although in the orthonormal case it is well understood and has elegant
theories~\cite{Donoho, AntFan, CandOver}, its shrinking and thresholding are not direct for a general regression matrix, and it suffers some problems in both selection and estimation~\cite{ZouHas, CandSta, Zhao}.
There has been a large and rapidly growing body of literature for the lasso studies over the past few years.
The efficient procedures proposed for solving the lasso include the well known LARS (Efron \emph{et al.}~\cite{Efron}), the homotopy method (Osborne \emph{et al.}~\cite{Osb}), and a recently re-discovered iterative algorithm (Fu~\cite{shooting}
Daubechies \emph{et al.}, Friedman \emph{et al.}~\cite{FHT}, Wu \& Lange~\cite{Lange}).
As for the theoretical aspects of the lasso, we refer to Knight \& Fu~\cite{Knight}, Zhao \& Yu~\cite{Zhao}, Donoho \emph{et al.}~\cite{DonohoSta}, Bunea \emph{et al.}~\cite{Bunea}, Zhang \& Huang~\cite{ZhangHuang}, etc. for  asymptotic and nonasymptotic results.
Various extensions and modifications to lasso have also been proposed, such as the grouped lasso (Yuan \& Lin~\cite{Yuan}), the Dantzig selector (Cand\`{e}s and Tao~\cite{Cand}), the adaptive lasso (Zou~\cite{Zou}), and the relaxed lasso (Meinshausen \& Yu~\cite{Relax}).

This paper aims to improve the na\"{\i}ve $l_1$-penalty, by using nonconvex penalties, to achieve an effective and efficient procedure for model selection and shrinkage.
The rest of the paper is organized as follows.
Section 2 provides a mechanism to  borrow the rich nonconvex penalties in the orthogonal design to solve the general problem.
From the point of view of thresholding rules, Section 3 constructs the thresholding-based iterative selection procedures (TISP) for model selection and successfully builds the convergence theorem.
Section 4 investigates the theoretical properties of the selection and the estimation via  TISPs nonasymptotically.
In Section 5, we carry out an empirical study of  TISP design which leads us to a novel Hybrid-TISP proposed  based on hard-thresholding and ridge-thresholding.
It provides a fusion between the $l_0$-penalty and the $l_2$-penalty, and adaptively achieves the right balance between shrinkage and selection in statistical modeling.
In practice, Hybrid-TISP shows superior performance in both test-error and sparsity. Section 6 gives a real data example.
All technical details are left to the Appendices.

\section{Motivation -- From Orthogonal Designs to Non-orthogonal Designs}
\label{secmotiv}
We consider the penalized regression problem
\begin{eqnarray}
\min_\bsbb \frac{1}{2} \|\bsbX \bsbb- \bsby\|_2^2 + P(\bsbb; \lambda)(\triangleq f(\bsbb)), \label{oriprob}
\end{eqnarray}
where $\bsbX=\left[\bsbx_1, \bsbx_2, \cdots, \bsbx_p\right]$ is the regression matrix, $\bsby\in R^n$ is the response vector, and $P(\bsbb; \lambda)$ represents the penalty with $\lambda$ as the regularization parameter. Here $p$ may be greater than $n$.
In this paper, we assume $\bsbb$ is sparse, and use \eqref{oriprob} for predictive learning. Although predictor error or accuracy is our first concern, we prefer to obtain a parsimonious model that is more interpretative in practice and  is consistent with Occam's razor. Usually $P$ is assumed to be an additive penalty in the sense that $P(\bsbb; \lambda)$ is obtained by a univariate $P$: $P(\bsbb; \lambda)=\sum P(\beta_i;\lambda)$. This sparsity problem has wide applications in variable selection, functional data analysis, graphical modeling, compressed sensing, and so on.

If $P(\bsbb; \lambda)=\lambda \|\bsbb\|_1$, then \eqref{oriprob} is the lasso~\cite{Tib}, a basic and popular method in variable selection. However, although the $l_1$-norm provides the best convex approximation to the $l_0$-norm and is computationally efficient, the lasso cannot handle collinearity~\cite{ZouHas} and may result in  inconsistent selection (cf. the irrepresentable conditions~\cite{Zhao}) and introduce extra bias in estimation~\cite{Relax}.

On the other hand, if we concentrate on orthogonal designs only, i.e., $\bsbX^T \bsbX=\bsbI$, like in wavelets, $l_1$ is far from the only choice. There are established theories and algorithms for various types of (nonconvex) penalties. 
\\

\noindent \textbf{Example 1. Hard-penalties}.
\begin{enumerate}
\item $P(\theta; \lambda)= \begin{cases} -\theta^2/2+\lambda |\theta|,  \mbox{ if } |\theta|<\lambda\\ \lambda^2/2,  \mbox{ if } |\theta|\geq \lambda \end{cases}$, due to Antoniadis~\cite{ant97}.
\item $P(\theta; \lambda)=\lambda^2/2 \cdot I_{\theta\neq 0}$, which is in fact the $l_0$-penalty.
\item $P(\theta;\lambda)=\begin{cases}\lambda |\theta|, \mbox{ if } |\theta|<\lambda\\ \lambda^2/2, \mbox{ if }|\theta|\geq \lambda \end{cases}$,  due to Fan~\cite{fan97}.
\end{enumerate}
It is interesting to note that all three lead to the same estimator obtained by hard-thresholding.

\noindent \textbf{Example 2. SCAD-penalty}.
${P'}(\theta;\lambda) = \begin{cases} \lambda, \mbox{ if } \theta\leq\lambda\\ (a\lambda-\theta)/(a-1), \mbox{ if } \lambda<\theta\leq a\lambda\\ 0, \mbox{ if } \theta > a \lambda\end{cases}$\\
for $\theta>0$ and $a>2$. The default choice of $a$ is $3.7$, based on a Bayesian argument (Fan~\cite{AntFan}). \\

\noindent \textbf{Example 3. Transformed $l_1$-penalty}.
$P(\theta;\lambda)={\lambda b |\theta|}/{(1+b|\theta|)}$ for some $b>0$, due to Geman \& Reynolds~\cite{transfL1}.\\

In this simplified setup, (a) the fitting part of the penalized regression \eqref{oriprob} is  \textbf{separable} in this case, which means we only need to deal with the univariate case, if $P$ is also separable (which is true in general);  (b) even if $P$ is \textbf{nonconvex}, it  still often results in a \emph{unique} solution.

One of our main goals in this paper is to borrow these rich results in the orthogonal design to help us solve the general problem \eqref{oriprob}. We use the following mechanism to achieve this.
Define
\begin{eqnarray}
g(\bsbb,\bsbg)=\frac{1}{2}\|\bsbX\bsbg-\bsby\|_2^2 + P(\bsbg;\lambda)+\frac{1}{2} <(\bsbI-\bsbSig)(\bsbg-\bsbb), \bsbg-\bsbb>.
\label{gdef}
\end{eqnarray}
Here $<\bsb{a}, \bsb{b}>=\bsb{a}^T \bsb{b}$, $\bsbSig=\bsbX^T\bsbX$.

Given $\bsbb$, minimizing $g$ over $\bsbg$ is equivalent to
\begin{eqnarray}
\arg \min_\bsbg \frac{1}{2} \left\|\bsbg - \left [(\bsbI-\bsbSig)\bsbb+\bsbX^T\bsby\right ]\right\|_2^2 + P(\bsbg; \lambda). \label{optovergamma}
\end{eqnarray}
In contrast to \eqref{oriprob}, this problem has an orthogonal design --- as mentioned earlier this is easier to handle both  in computation and in theory. For example, we may adopt some nonconvex penalties, and they still result in a unique solution of $\bsbg$.

Given $\bsbg$, minimizing $g$ over $\bsbb$ is equivalent to
\begin{eqnarray}
\arg \min_\bsbb \frac{1}{2} <(\bsbI - \bsbSig)\bsbb, \bsbb-2\bsbg>. \label{optoverbeta}
\end{eqnarray}
Taking its derivative with respect to $\bsbb$ gives $(\bsbI-\bsbSig)(\bsbb-\bsbg)=\bsb{0}$, from which it follows that $\bsbb=\bsbg$ if $\|\bsbSig\|_2<1$.
Note that \eqref{optoverbeta} is a convex optimization.  Therefore, the optimal value of $g$ is always achieved at $\bsbg=\bsbb$ if $\bsbX$ is scaled down properly.

The connection to the original problem is now clear: it is easy to verify $\min_\bsbb g(\bsbb, \bsbb)$ is equivalent to $\min_\bsbb f(\bsbb)$. The advantage of optimizing $g$ instead of $f$ is that given $\bsbb$, the problem is orthogonal and separable in  $\bsbg$, and we can adopt far more flexible penalties in the algorithm design, including the nonconvex ones.

\section{Thresholding-based Iterative Selection Procedures (TISP)}
\subsection{Thresholding Rules and Penalties}
\label{subsecthpconst}
As the title suggests, our starting point in this paper is  thresholding rules rather than different forms of the penalty function. One direct reason is that different $P$'s may result in the same estimator and the same thresholding, say, in the situation of hard-thresholding~\cite{ant97, fan97}.
Moreover, starting with thresholding functions facilitates the computation (as will be shown in the next subsection).
Besides, there is also a universal connection between  thresholding rules and  penalty functions that we will investigate in this subsection. For convenience, we consider the univariate case only.

A thresholding function, denoted by $\Theta(\cdot; \lambda)$, with $\lambda$ as a parameter, is required to satisfy:
\begin{enumerate}
\item $\Theta(\cdot;\lambda)$ is an odd function. ($\Theta_+(\cdot;\lambda)$ is used to denote the $\Theta(\cdot;\lambda)$ restricted to $R_+=[0,\infty)$.)
\item $\Theta$ is a shrinkage rule: $0\leq \Theta_+(t;\lambda)\leq t, \forall t\in R_+$.
\item $\Theta_+$ is nondecreasing on $R_+$, and $\Theta_+(t;\lambda)\rightarrow\infty$ as $t\rightarrow\infty$.
\end{enumerate}
In addition, it is natural to have $\Theta_+(t;\lambda)=0, 0\leq t \leq \tau$ for some $\tau \geq 0$.


Given a thresholding rule $\Theta(\cdot; \lambda)$, a penalty function can be obtained from the following three-step construction. First, define
$$
\Theta^{-1}(u;\lambda)=\sup\{t:\Theta(t;\lambda)\leq u\} \mbox{ and } \Theta^{-1}(-u;\lambda)=-\Theta^{-1}(u;\lambda),
$$
for any $u\in R_+$. Then define
\begin{equation}
s(u;\lambda)\triangleq \Theta^{-1}(u;\lambda)-u, \forall u. \label{defofs}
\end{equation}
Finally, let $P$ be a continuous and positive penalty defined by
\begin{eqnarray}
P(\theta;\lambda)=\int_0^{|\theta|} s(u;\lambda) du. \label{constrP}
\end{eqnarray}
Antoniadis~\cite{antrev} showed the following result for this constructed $P$.
\begin{proposition}
\label{uniqsol}
The minimization problem $\min_\theta  (t-\theta)^2/2 + P(\theta;\lambda)$ has a unique optimal solution $\hat\theta=\Theta(t;\lambda)$ for every $t$ at which $\Theta(\cdot;\lambda)$ is continuous.
\end{proposition}

In addition, if we define
$
\psi(t)=t-\Theta(t), 
$ 
then it is the psi-function for defining M-estimators; see~\cite{antrev, gannaz}.


Note that \eqref{constrP} is not the only way to construct a penalty that leads to $\Theta$ in solving the optimization. For example, in the situation of hard-thresholding, in addition to the continuous penalty
\begin{eqnarray}
P=\lambda^2/2-(|\theta|-\lambda)^2 1_{|\theta|<\lambda}/2 \label{hardPcont}
\end{eqnarray} constructed via \eqref{constrP},
\begin{eqnarray}
P(\theta;\lambda)=\begin{cases}\lambda |\theta|, \mbox{ if } |\theta|<\lambda\\ \lambda^2/2, \mbox{ if }|\theta|\geq \lambda \end{cases}, \mbox{ and }  P(\theta;\lambda)=\frac{\lambda^2}{2}\cdot 1_{\theta\neq 0} \label{hardPs}
\end{eqnarray}
are also valid choices~\cite{ant97, fan97}. In some sense, \eqref{hardPcont} may be considered as a continuous version of the discrete $l_0$-penalty.

\subsection{TISP and Its Convergence}
Now we go back to the mechanism introduced in Section \ref{secmotiv} for the  penalized multivariate regression problem \eqref{oriprob}, with $P$ constructed from a given  thresholding function $\Theta$.
Solving \eqref{optovergamma} yields $
\bsbg = \Theta((\bsbI-\bsbSig)\bsbb+\bsbX^T\bsby;\lambda)
$.
Seen from \eqref{optoverbeta}, our iterates simplify to
\begin{eqnarray}
\bsbb^{(j+1)} = \Theta((\bsbI-\bsbSig)\bsbb^{(j)}+\bsbX^T\bsby;\lambda).
\label{tisp}
\end{eqnarray}
This iterative procedure  is referred to as the \textbf{Thresholding-based Iterative Selection Procedure} (\textbf{TISP}).
TISP provides a feasible way to tackle the original optimization \eqref{oriprob}. It is a simple procedure that does not involve any complicated operations like matrix inversion.

There are rich examples for the procedure defined by \eqref{tisp}. (a) Using a soft-thresholding in \eqref{tisp}, we immediately obtain the iterative algorithm (in vector form) for solving the lasso problem where $P(\bsbb; \lambda)=\lambda \|\bsbb\|_1$~\cite{Daub}. In fact, the asynchronous updating of \eqref{tisp}  leads exactly to the component-by-component iteration referred to as the coordinate decent algorithm (see Friedman \emph{et al.}~\cite{FHT}). The corresponding pathwise algorithm has been considered to be the fastest in solving the lasso problem to date, especially when $p>n$.
(b) If we substitute hard-thresholding for $\Theta$, seen from \eqref{hardPs}, it is an alterative optimization for solving the  penalized regression  with
 $$P=c\cdot \sum_i 1_{\beta_i\neq 0}=c\cdot \|\bsbb\|_0,$$ i.e., the $l_0$-penalized regression problem.
(c) We can also replace the hard-thresholding by the more smoothed SCAD to reduce instability.
(d) Finally, it is worth mentioning that TISP may also include the ridge penalty
$P(\bsbb;\lambda)=\lambda \|\bsbb\|_2^2/2$, if we set
\begin{eqnarray}
\Theta(t;\lambda)=\frac{t}{1+\lambda},\label{ridgethfunc}
\end{eqnarray} thanks to the generic definition of a thresholding function.\\

Obviously, if $\bsbSig$ is nonsingular, and so $n>p$, the TISP mapping is  a contraction and thus the sequence $\bsbb^{(j)}$ converges to a stationary point of \eqref{oriprob}. We would like to apply TISP to large $p$ problems as well where $\bsbSig$ is singular --- a surprising fact  is, however, that
TISP  may not be a \emph{nonexpansive} operator\footnote{An operator $T$ is called nonexpansive~\cite{Browder} if $\|T(x)-T(y)\|\leq \|x-y\|$ for \emph{any} $x,y$. Obviously, the hard-thresholding function is not nonexpansive.} for most thresholdings (except soft-thresholding), let alone a contraction. The following studies cover the large $p$ case ($p>n$).  We use $\mu(\bsb{A})$ to represent an arbitrary  singular value of matrix $\bsb{A}$, and $\mu_{\max}(\bsb{A})$ ($\mu_{\min}(\bsb{A})$)  the max (min) of $\mu(\bsb{A})$, respectively. 

Without loss of generality, suppose the penalty function  defined by \eqref{constrP} satisfies the bounded curvature condition (BCC) for some symmetric matrix $\bsbH$:
\begin{eqnarray}
P(\bsbb+\bsb{\Delta}; \lambda)\geq P(\bsbb; \lambda)+<\bsb{\Delta},\bsb{s}>-\frac{1}{2}\bsb{\Delta}^T\bsbH\bsb{\Delta}, \forall \bsb{\Delta} \in R^p \label{bcc}
\end{eqnarray}
where $\bsb{s}=\bsb{s}(\bsbb;\lambda)$ is given by \eqref{defofs}.
Many  thresholding rules of practical interest including Example 1-3 satisfy the BCC with a positive semi-definite $\bsbH$. For instance, for soft-thresholding, $\bsbH=\bsb{0}$ since $\|\bsbb\|_1$ is convex; for hard-thresholding, $\bsb{H}=\bsb{I}$; for SCAD-thresholding, we can take $\bsbH=\bsbI/(a-1)$ (recall that the parameter $a$ is assumed to be greater than $2$ in Example 2, and so $\bsbH$ is positive definite).

\begin{theorem}
\label{conv}
Given the TISP \eqref{tisp}, if $\mu_{\max}(\bsbSig)\leq 1\vee (2-\mu_{\max}(\bsbH))$, then
\begin{eqnarray}
f(\bsbb^{(j)}) \geq f(\bsbb^{(j+1)}). \label{optfval}
\end{eqnarray}
Moreover, if $\mu_{\max}(\bsbSig)< 1\vee (2-\mu_{\max}(\bsbH))$, there exists a constant $C>0$, dependent on $\bsbX$, $\bsbH$ only, such that
\begin{eqnarray}
f(\bsbb^{(j)})- f(\bsbb^{(j+1)})\geq C\cdot \|\bsbb^{(j)}-\bsbb^{(j+1)}\|_2^2. \label{asympreg}
\end{eqnarray}
\end{theorem}

Therefore, for an arbitrary $\bsbX$, we can use  TISP of the following form in practice
\begin{eqnarray}
\bsbb^{(j+1)} = \Theta\left(\left(\bsbI-\frac{1}{k_0^2}\bsbSig\right)\bsbb^{(j)}+\frac{1}{k_0^2}\bsbX^T\bsby;\frac{\lambda}{k_0^2}\right),
\label{tispvar}
\end{eqnarray}
where $k_0=\mu_{\max} (\bsbX)=\|\bsbX\|_2$, although larger values of $k_0$ generally lead to faster convergence.
Applying Theorem \ref{conv} to some interesting special cases gives the following corollaries.

\begin{corollary}
\label{softconv}
Suppose $\Theta$ is soft-thresholding. If $\mu_{\max}(\bsbX) < \sqrt 2$, then \eqref{asympreg} holds.
\end{corollary}

\begin{corollary}
\label{hardconv}
Suppose $\Theta$ is hard-thresholding. If $\mu_{\max}(\bsbX) \leq 1$, then \eqref{optfval} holds; further, if  $\mu_{\max}(\bsbX) < 1$, then \eqref{asympreg} is true.
\end{corollary}

\begin{corollary}
\label{scadconv}
Suppose $\Theta$ is SCAD-thresholding. If $\mu_{\max}(\bsbX) < \sqrt{2-\frac{1}{a-1}}$, then \eqref{asympreg} holds.
\end{corollary}

Corollary \ref{softconv} generalizes the lasso result by Daubechies \emph{et al.}~\cite{Daub}, and coincides with our previous study~\cite{Shethesis}. Corollary \ref{scadconv} covers the orthogonal case, since SCAD assumes $a>2$  and thus $\sqrt{2-\frac{1}{a-1}}>1$.
Finally, it is worth pointing out that TISP may not always be an MM algorithm~\cite{HuntLang} like the  LLA method by Zou \& Li~\cite{ZouLi}. Take the SCAD-thresholding as an example: when $1< \|\bsbX\|_2 < \sqrt{2-\frac{1}{a-1}}$, $g$ defined by \eqref{gdef} does not majorize  $f$ but TISP converges.  Theorem \ref{conv} implies that if $\bsbX$ is scaled down properly (which does not affect the variable selection), $f(\bsbb^{(j)})$ is nonincreasing all the time during the iteration process.

We can easily show a result similar to Zou \& Li~\cite{ZouLi}:
\begin{proposition}
\label{convstat}
Suppose $\mu_{\max}(\bsbSig)< 1\vee (2-\mu_{\max}(\bsbH))$.
Give an initial point $\bsbb{(0)}$, if $\bsbb^*$ is a limit point of the TISP sequence $\bsbb^{(j)}$, then $\bsbb^*$ is a stationary point of $f(\bsbb)$ \eqref{oriprob}, or equivalently, a fixed point of \eqref{tisp}.
\end{proposition}

Denote by $F$ the set of the fixed points of TISP. That is, given any $\bsbb^*\in F$, it satisfies the  implicit equation
\begin{eqnarray}
\bsbb = \Theta((\bsbI-\bsbSig)\bsbb+\bsbX^T\bsby;\lambda),
\label{tispkkt}
\end{eqnarray}
referred to as the \emph{$\Theta$-equation}.
Clearly, local minima of $f$ are fixed points of \eqref{tispkkt}.  In the next section, we will perform an nonasymptotic study of the good properties of the points in $F$. Here, we give the following optimality result.
\begin{proposition}
\label{optimality}
Let $\bsbb^*\in F$ and suppose $\mu_{\max}(\bsbH)\leq 1$. If $\mu_{\max}(\bsbH)\leq \mu(\bsbSig)\leq 2-\mu_{\max}(\bsbH)$, then $\bsbb^*$ is a global minimizer of $f$.
\end{proposition}

Although the fact that nonconvex penalties often result in a unique optimal solution in the \emph{orthogonal} design is well known, this proposition states (novelly) that the same conclusion holds as long as $\bsbX$ is not too far  from orthogonal (characterized in terms of $\bsbH$). For instance, for SCAD thresholding and penalty, TISP necessarily leads to the global minimum of $f$, provided $\frac{1}{\sqrt{a-1}}\leq \mu(\bsbX)\leq \sqrt{2-\frac{1}{a-1}}$, or $0.61\leq \mu(\bsbX)\leq 1.27$ when $a=3.7$ (the default choice in SCAD--see Example 2), given any initial point $\bsbb^{(0)}$. In summary, TISP is a successful algorithm for solving the penalized regressions for a general design matrix.

\subsection{Related Work}
\label{subsecrelwork}
The main contribution of this paper is to consider a new class of \emph{$\Theta$-estimators}  defined by the $\Theta$-equation \eqref{tispkkt} for model selection and shrinkage, which
can be naturally computed by TISP, and are associated with penalized regressions---in particular, the penalty $P$ can be constructed via the three-step procedure introduced in Section \ref{subsecthpconst}.
More generally, the $\lambda$ in \eqref{tispkkt} can be component-specific.  For example, if $\bsbX$ is not column-normalized, we may use
\begin{eqnarray}
\bsbb = \Theta\left(\left(\bsbI-\bsbSig\right)\bsbb+\bsbX^T\bsby;\bsb{\lambda}\right),  \label{tispkktvec}
\end{eqnarray}
where $\bsb{\lambda}=\left[\begin{array}{cccc} \lambda \|\bsbx_1\|_2 &  \lambda \|\bsbx_2\|_2 &\cdots &  \lambda \|\bsbx_p\|_2\end{array} \right]^T$ and $\lambda$ is a regularization parameter.
With a carefully designed $\Theta$, we obtain a good estimator with both accuracy and sparsity, as will be shown in Section \ref{subsechybridtisp}.  The corresponding penalty is, not surprisingly, nonconvex, which indicates the difficulty of this NP-hard problem.

Nonconvex penalties have been successfully used in real-world applications like high-dimensional nonparametric modeling~\cite{AntFan}, survival analysis~\cite{Caisurv}, and  microarray data analysis~\cite{scadreg, scadsvm}, where they achieve outstanding performance. The numerical optimization has been a challenging and intriguing  problem.
In the context of wavelet denoising where $\bsbX\bsbX^T=\bsbI$, Antoniadis \& Fan~\cite{AntFan} proposed the ROSE to approximately solve the minimization problem for a wide class of nonconvex penalties. They also introduced the graduated nonconvexity (GNC) algorithm, developed in image processing; it has a number of tuning parameters and is computationally intensive.
Fan and Li~\cite{FanLi} then proposed a generic local quadratic approximation (LQA) algorithm by solving a series of $l_2$-penalized problems. Like ridge regression, this approach does not intrinsically yield  zeros, and setting a small cutoff value during iteration has been shown to be too greedy. A refined version is the perturbed LQA suggested by Hunter \& Li~\cite{HunLi} to avoid numerical instability.   The perturbation parameter needs to be chosen very carefully in implementation since it affects the sparsity of the solution as well as the speed of convergence. Recently, Zou \& Li~\cite{ZouLi} proposed a new local linear approximation (LLA) which significantly improves the LQA. Explicit sparsity is attained by solving a weighted lasso problem at each iteration step. (Note that our TISP does a simple thresholding at each step.) One-step SCAD estimator is advocated. Our empirical studies show that this one-step convex approximation has limited power in finite samples. Although the estimate is sparser than using the plain $l_1$-penalty, it may result in misleading models with poor prediction error. See Section \ref{secdesign} for detail.

Using thresholding rules to define $\Theta$-estimators shares similarities to the studies of $M$-estimators of  $\psi$-type in robust regression. Most $M$-estimators were proposed in the form of $\psi$-functions but not based on  loss functions, such as  Huber's, Hampel's three-part, and Tukey's bisquare $M$-estimators. Indeed, we find an interesting connection between these two fields.
Assume a mean shift outlier model, $\bsby = \bsbX\bsbb + \bsbg +\bsbeps, \bsbeps\sim N(0, \sigma^2 \bsbI)$,  where $n>p$ and $\bsbg\in R^{n}$ is sparse. If $\bsbg_i$ is nonzero, case $i$ is an outlier.
Let $\bsbH=\bsbX(\bsbX^T\bsbX)^{-1}\bsbX^T$ be the hat matrix and suppose its spectral decomposition is given by $\bsbH=\bsbU \bsbD \bsbU^T$. Define an index set $c=\{i: D_{ii} = 0\}$ and $\bsbU_c$ is formed by taking the corresponding columns of $\bsbU$. Then a reduced model can be obtained from the mean shift outlier model
\begin{eqnarray}
\tilde\bsby=\bsbA\bsbg +  {\bsbeps'}, \bsbeps'\sim N(\bsb{0}, \sigma^2\bsbI_{(n-p)\times (n-p)}), \label{outlSpa}
\end{eqnarray}
where    $\tilde\bsby=\bsbU_c^T\bsby$, $\bsbA = \bsbU_c^T $ (of dimension ${(n-p)\times n}$). Applying $\Theta$-TISP to this `large-$p$' sparsity problem gives the iteration $\bsbg^{(j+1)} = \Theta(\bsbA^T \tilde\bsby/k_0^2  + (\bsbI-\bsbA^T\bsbA/k_0^2)\bsbg^{(j)}; \bsb{\lambda}/k_0^2)$ with  $k_0 = \mu_{max}(\bsbH)$, or
 $$\bsbg^{(j+1)} = \Theta(\bsbH \bsbg^{(j)} + (\bsbI-\bsbH)\bsby; \bsb{\lambda}).$$
After getting $\hat\bsbg$ from TISP, we can estimate $\bsbb$ by $\hat\bsbb = (\bsbX^T\bsbX)^{-1}\bsbX^T(\bsby-\hat\bsbg)$.
Simple algebra shows that this special $\Theta$-TISP solves an $M$-estimation problem associated with $\psi$, if ($\Theta$, $\psi$) satisfies $\Theta(t; \lambda) + \psi(t; \lambda)=t.$
It is well known that Huber's method (or equivalently, Soft-TISP, which corresponds to using a convex $l_1$-penalty on $\bsbg$) behaves poorly in outlier detection even for moderate leverage points.  Instead, redescending $\psi$-functions are advocated, which corresponds to using nonconvex penalties for the sparsity problem of \eqref{outlSpa}.

\section{Selection and Estimation via TISP}
\label{sectheo}
TISP provides a very simple way to do variable selection via penalized regressions.
In this section, we will perform a theoretical study of the variable selection and coefficient estimation by TISPs based on different thresholdings. Our results are nonasymptotic.
\subsection{Assumptions on $\Theta$}
Given $\Theta(\cdot;\lambda)$, denote its thresholding value by $\tau(\lambda)$, i.e., $\Theta(t;\lambda)=0 \mbox{ } \forall t: |t| < \tau$ and $\Theta(t;\lambda)\neq 0  \mbox{ for } |t| > \tau$. For example, $\tau(\lambda)=\lambda$ in soft-, hard-, and SCAD-thresholdings, but is not so for the transformed $l_1$. Assume $\tau>0$. To ease our TISP study based on the $\Theta$-equation \eqref{tispkkt}, we define another version of $\bsbs$, called the generalized sign. Introduce
$$
\gSgn(u;\lambda)=\{s\in R: \Theta(u+\tau s;\lambda)=u\}\quad \mbox{ if }  u\in \mbox{ran}(\Theta),
$$
and $\gSgn(u;\lambda)=\{0\}$ otherwise, where $\mbox{ran}(\Theta)$ is the range of $\Theta$; $\gsgn(u;\lambda)$ is  used to denote a specific element in $\gSgn(u;\lambda)$. The vector versions of $\gSgn$ and $\gsgn$ can be defined correspondingly. Clearly if $u = \Theta(t;\lambda)$ then $t=u+\tau \gsgn(u;\lambda)$ for some $\gsgn(u;\lambda)\in \gSgn(u;\lambda)$.

As a demonstration, if $\Theta(\cdot;\lambda)$ is soft-thresholding, $\tau=\lambda$ and $\widetilde{\mbox{Sgn}}(\bsbb)=\{\bsb{s}: s_i=1 {\mbox{ if }} \beta_i>0, s_i=-1 {\mbox{ if }} \beta_i<0, {\mbox{ and }} s_i\in[-1,1] {\mbox{ if }} \beta_i=0\}$. Thus
now $\gSgn(\bsbb)$  is the subdifferential of $\|\bsbb\|_1$, and $\gsgn(\bsbb)$ is a subgradient~\cite{Shim}.
For hard-thresholding, $\widetilde{\mbox{Sgn}}(\bsbb)=\{\bsb{s}: s_i=0 {\mbox{ if }} \beta_i\neq 0,  s_i\in[-1,1] {\mbox{ if }} \beta_i=0\}$. $\widetilde{\mbox{Sgn}}$ and $\widetilde{\mbox{sgn}}$ are called generalized signs due to the following fact. 

\begin{proposition}
\label{gsgn}
Suppose $\Theta(\cdot;\lambda)$ is sandwiched by soft- and hard-thresholdings, $\Theta_S(\cdot;\tau)$ and $\Theta_H(\cdot;\tau)$, i.e.,
\begin{eqnarray}
(\Theta_S)_+(t;\tau)\leq \Theta_+(t;\lambda) \leq (\Theta_H)_+(t;\tau), \forall t \in R_+. \label{sandwichcond}
\end{eqnarray}
Then $0\leq \gsgn(u)\leq 1$ if $u>0$, $-1\leq \gsgn(u)\leq 0$ if $u<0$, and $\gsgn(0)\in [-1, 1]$.
\end{proposition}

This proposition is easy to prove from the non-decreasing property of $\Theta$. Throughout the rest of the section, we assume $\Theta$ always satisfies the sandwiching condition \eqref{sandwichcond}.
By the definition of the generalized signs, \eqref{tispkkt} is equivalent to
$
\bsbSig \bsbb=\bsbX^T \bsby - \tau\gsgn(\bsbb;\lambda), 
$ 
for some $\gsgn(\bsbb;\lambda)\in \gSgn(\bsbb;\lambda)$.
We study the TISP estimate based on the scaled form \eqref{tispvar}.
Let $\hat\bsbb$ be a fixed point of \eqref{tispvar} and suppose $\tau(\lambda)=c\tau(\lambda/c)$ for any $c>0$. Then  the $\Theta$-equation for this TISP estimate can be rewritten as
\begin{eqnarray}
\bsbSig \hat \bsbb=\bsbX^T \bsby - \tau\gsgn(\hat \bsbb;\lambda/k_0^2), \label{kktgsgn}
\end{eqnarray}
where $k_0=\|\bsbX\|_2$.

\subsection{Sparsity Recovery}
Recall that $\bsby = \bsbX \bsbb + \bsbeps$, $\bsbeps \sim N(\bsb{0}, \sigma^2 \bsb{I})$, and $\bsbb$ is sparse. Let $z=\{i: \beta_i = 0\}$,  $nz=\{i: \beta_i \neq 0\}$, $d_z=|z|$, $d_{nz}=|nz|$. To study the sign-consistency of a TISP estimate, we denote by $p_s$ the probability of successful sign recovery, that is, the probability that there exists a $\hat\bsbb\in F$ such that ${\mbox{sgn}}(\hat\bsbb)={\mbox{sgn}}(\bsbb)$.

To simplify asymptotic discussions, we assume $\bsbX$ has been scaled to have all column $l_2$-norms equal to $\sqrt n$. Define $\bsbSig^{(s)}=\bsbSig/n$.
To get a better form of the bounds for $p_s$, we define two quantities $\mu = \mu_{\min}(\bsbSig_{nz,nz}^{(s)})$ and $\kappa \triangleq \underset{ i \in z}{\max} \| \bsbSig_{i,nz}^{(s)} \|_2 / \sqrt{ d_{nz}}$. Intuitively, $\kappa$ measures the `mean' correlations between the relevant predictors and the irrelevant predictors. The following nonasymptotic result is always true regarding the  selection via TISP.

\begin{theorem}
\label{selbytisp}
Assume $\mu \geq  \kappa d_{nz}$, $\mu>0$ and $\min|\bsbb_{nz}| \geq \frac{d_{nz}\tau }{ n\mu}$, then
\begin{eqnarray}
p_s \geq [1-2\Phi(-M)]^{d_z} [1-2\Phi(-L)]^{d_{nz}}, \label{selbndtisp}
\end{eqnarray}
where $M=\left(1- \frac{\kappa   d_{nz}}{\mu} \right)\frac{\tau}{\sqrt n\sigma}$, $L
= \frac{\sqrt{\mu n}}{\sigma} \left(   \min|\bsbb_{nz}| - \frac{\tau d_{nz}}{ {\mu n}}\right)$, and $\Phi$ is the standard normal distribution.

\end{theorem}

\begin{corollary}
\label{selbytispcor}
Under the conditions of Theorem \ref{selbytisp}, we have
\begin{eqnarray}
 1 - p_s \leq 2d_z\varphi(M)/M +2 d_{nz}\varphi(L)/L,\label{selbndtisp2}
\end{eqnarray}
where $\varphi$ is the standard normal density.
\end{corollary}

Clearly, the size of $\kappa$ is very important.
A small value of $\kappa$  {weakens} the interference of $\bsbX_z$ and $\bsbX_{nz}$ and helps recover the sparsity correctly.
We can also use this theorem to explore some asymptotics. (i) Assume $\bsbb$, $d_z$, and $d_{nz}$ are fixed, $n \rightarrow \infty$, then  under some
regularity conditions we get: if $\tau /\sqrt n \rightarrow \infty$ and $\tau  / n
\rightarrow 0$, then  TISP is sign consistent. This result in the Soft-TISP (lasso) case coincides with other studies
like~\cite{Knight,Zhao}.
(ii) Suppose $\bsbb_{nz}$ and $d_{nz}$ are fixed,
$n, d_z \rightarrow \infty$, and $\mu  \geq  (1+\epsilon)\kappa d_{nz}$ for some $\epsilon>0$.
Then TISP can successfully recover the sparsity pattern of $\bsbb$ if
$d_z\varphi(M)/M \rightarrow 0$ and $\tau/n \rightarrow 0$, which only requires $n$
to grow faster than $\log d_{z}$.

Unfortunately, the regularity condition $\mu \geq  \kappa d_{nz}$ cannot be removed in general. In the lasso case, it is a version of the irrepresentable conditions~\cite{Zhao}. (We took this more restrictive form because it is more intuitive and leads to more nice-looking bounds in \eqref{selbndtisp} and \eqref{selbndtisp2}.) However, for hard-thresholding-like $\Theta$'s, this is unnecessary and we can obtain stronger results.

We say that $\Theta$ belongs to the \emph{hard-thresholding family} if \begin{eqnarray}
\Theta(t;\lambda) =t, \forall t: |t|>c\cdot \tau, \label{hardcond}
\end{eqnarray}
for some constant $c\geq 1$. Hard-thresholding and SCAD-thresholding  are two examples with $c=1$, $a$ respectively. Unlike soft-thresholding, they do not introduce bias for large nonzero components.

\begin{theorem}
\label{selbyhard}
Suppose $\Theta$ belongs to the hard-thresholding family and $\min|\bsbb_{nz}| \geq c\tau/k_0^2$.
Then
\begin{eqnarray}
p_s \geq [1-2\Phi(-M')]^{d_z} [1-2\Phi(-L')]^{d_{nz}}, \label{selbndhard}
\end{eqnarray}
where $M'=\frac{c\tau}{\sqrt n \sigma }$, $L'= \frac{\sqrt{\mu n} }{\sigma} \left(  \min|\bsbb_{nz}| - \frac{c\tau}{ k_0^2}\right)$.
\end{theorem}

\begin{corollary}
\label{selbyhardcor}
Under the conditions of Theorem \ref{selbyhard}, we have
\begin{eqnarray}
1 - p_s \leq 2 d_z\varphi(M')/M' +2 d_{nz}\varphi(L')/L'. \label{selbndhard2}
\end{eqnarray}
\end{corollary}

\eqref{selbndhard} is strictly better than the bound in \eqref{selbndtisp} if $c<d_{nz}k_0^2/(\mu n)$, or $c<d_{nz} {\mu_{\max}(\bsbSig^{(s)})}/{\mu_{\min}(\bsbSig^{(s)}_{nz})}$, which is usually true for both hard- and scad-thresholding.
Therefore the TISP induced by a $\Theta$ in the  hard-thresholding family can achieve better performance in variable selection.\footnote{Note that, however, the regularization parameters are generally tuned to reduce the test error.}
This  will be verified empirically in the next section.
Note that although in the orthogonal case, hard-thresholding and soft-thresholding give exactly the same zeros, they result in very different sparsity patterns  in our iterative procedure for a nonorthogonal $\bsbX$.

\subsection{Estimation Risk}
We obtain the following TISP risk bounds for any thresholding $\Theta$.
\begin{theorem}
\label{estbytisp}
Let $\nu = \mu_{\min} (\bsbSig_{z,z}^{(s)})$ and $\hat \bsbb\in F$. Define $R_{nz} = E(\|\bsbb_{nz} - \hat \bsbb_{nz}\|_2^2)$, and $R_{z} = E( \|\hat \bsbb_{z}\|_2^2)$. Suppose $\bsbSig$ is nonsingular. Then
\begin{eqnarray}
R_{nz} \leq \frac{3}{n}\left[
\frac{d_{nz}}{\mu}\sigma^2 + \frac{d_{nz}}{\mu^2}\frac{\tau^2}{n}  + \kappa^2  \frac{d_{z}
d_{nz}}{\mu^2}\cdot  n R_{z}\right].
\label{nzriskbnd}
\end{eqnarray}
And
\begin{eqnarray}
R_z \leq \frac{\sigma^2}{n}\frac{d_z^2 }
{\nu^2}  ( K_1 M + K_2 \frac{1}{M}) \varphi(M),
\label{zriskbnd}
\end{eqnarray}
where $M$ is defined as in Theorem \ref{selbytisp}, $K_1=6  \frac{1+\frac{1+{\kappa^2 d_{nz}^2}/{\mu^2} }{\left(1-{\kappa
d_{nz}}/{\mu}\right)^2} }{\left( 1 - \kappa^2 \frac{ d_z d_{nz} }{  \mu\nu }\right)^2}$,
$K_2=6 {\left( 1 - \kappa^2 \frac{ d_z d_{nz} }{  \mu \nu}\right)^{-2}}$ in which we assume $\kappa^2  \leq \frac{\mu\nu  }{  d_z d_{nz} } $ and $\mu\geq \kappa d_{nz}$.
\end{theorem}
This general result holds for any TISP estimate.
It is not difficult to show that in the previous setting of (i) where $\bsbb$ is fixed, we also obtain $R_z \rightarrow 0$ and $R_{nz} \rightarrow 0$ by the theorem if $\tau/\sqrt{n}\rightarrow \infty$ and $\tau/n\rightarrow \infty$.
Besides, under the  conditions stated in the theorem,
$M=\sqrt{2\log \frac{d_z^2}{n}+(1+\epsilon)\log\log \frac{d_z^2}{n}}$  is sufficient to ensure $R_z \rightarrow 0$ for any $\epsilon > 0$; since  $M \sim \sqrt{2 \log \frac{d_z^2}{n}} \leq \sqrt{2 \log d_z } $, Donoho \& Johnstone's classical work~\cite{Donoho} in the orthogonal design implies this risk bound can not be improved significantly in general.
We leave the TISP design problem  to the next section using an empirical study.

In the orthogonal case, we can show the oracle inequalities~\cite{Donoho} hold.
\begin{theorem}
\label{orthoestbytisp}
Suppose $\Theta$ satisfies the sandwiching condition \eqref{sandwichcond} and $\bsbX^T\bsbX=\bsbI$. Then
\begin{eqnarray}
E \| \hat\bsbb - \bsbb \|_2^2 \leq (1+\tau^2) \sum_1^n \min \left(
\frac{2\varphi(\tau)}{\tau}\sigma^2+\beta_i^2,\sigma^2 \right)
\end{eqnarray}
for any $\tau>1$. Consequently, when $\tau=\sqrt{2 \log n}$,
\begin{eqnarray}
E \| \hat\bsbb - \bsbb \|_2^2 \leq (2 \log n + 1) \left( \frac{\sigma^2}{\sqrt{\pi \log n}}+\sum
\min(\beta_i^2,\sigma^2)\right)
\end{eqnarray}
for any $n\geq 2$.
\end{theorem}

This nonasymptotic result covers  soft-,  hard-, and SCAD-thresholdings. It coincides with the classical soft-thresholding studies~\cite{Donoho} and is sharper than~\cite{AntFan, Zou}.  (A correction of Zou's oracle bound~\cite{Zou} is also given at the end of the proof; see Appendix \ref{appproofOrtho}.)

\section{TISP Designs: An Empirical Study}
\label{secdesign}
\subsection{A numerical study of TISPs}
\label{subsecdesign}
In this section, we demonstrate the empirical performance of TISPs by some simulation data.
Although there are rich choices about $\Theta$ in \eqref{tisp}, we focus on three basic TISPs only in this subsection.
In addition to the Soft-TISP, i.e., the lasso, we implemented Hard-TISP and SCAD-TISP, the thresholdings of which belong to the \emph{hard-thresholding family}. The parameter $a$ in SCAD-thresholding takes the default value, $3.7$, based on a Bayesian argument~\cite{FanLi}. As seen from the theoretical studies in   Section \ref{sectheo}, the last two should perform better than the lasso in variable selection.
In generating the solution path for a grid of $\lambda$-values,
we always set the initial point, $\bsbb^{(0)}$, to be  zero in Hard- or SCAD-TISP.
A natural search range for $\lambda$, seen from \eqref{tispkkt}, is $[0, \bsbX^T\bsby]$, if $\bsbX$ has been column normalized.
(Note that a pathwise algorithm with warm start, which takes the previous estimate associated with the old value of $\lambda$ as the initial point of the procedure
for the current value of $\lambda$, may be inappropriate for TISPs when nonconvex penalties are used. In fact, the solution path associated with a nonconvex penalty is generally not continuous in $\lambda$ and warm-start leads to bad solutions because of multiple local minima effects.)

For comparison, the one-step LLA method,  proposed by Zou \& Li~\cite{ZouLi} for  penalized likelihood models, is also included in our tests.
They showed good asymptotics about one-step SCAD when $n\rightarrow\infty$ and $p$ is fixed, and demonstrated its  performance in
various numerical examples.
The one-step LLA is actually a {weighted lasso} with weights constructed from the OLS estimate using different penalty functions.
According to our general result of weights in sparse regression~\cite{Shethesis}, it can achieve better sign consistency than the lasso as $n$ grows to infinity.
We are greatly interested in drawing a comparison between TISP and LLA since TISP also successfully solves the  penalized
regression problems.

We did experiments on two simulation datasets. Each dataset contains training data, validation data, and
test data.  We use $\#=$``$\cdot / \cdot / \cdot$" to denote the number of
observations in the training data, validation data, and test data. Let $\bsbSig$ be the correlation
matrix in generating $\bsbX$, i.e., each row of $\bsbX$ is independently drawn from $N(\bsb{0}, \bsbSig)$. We use
$(\{a_1\}^{n_1}, \cdots, \{a_k\}^{n_k})$ to denote the column vector made by $n_1$ $a_1$'s,
$\cdots$, $n_k$ $a_k$'s
consecutively in the following examples.\\
\textbf{Example 1.} $\#=20/100/200$, $d=8$, $\bsbb = (\{3\}^1, \{1.5\}^1, \{0\}^2, \{2\}^1,
\{0\}^3)$, $\Sigma_{ij}=\rho^{|i-j|}$ with $\rho=0.5$, $\sigma=2, 3, 5, 8$; the corresponding signal-to-noise variance ratio ($\bsbb^T\bsbSig\bsbb/\sigma^2$) is $5.31$, $2.36$, $0.85$, and $0.33$, respectively.\\
\textbf{Example 2.} $\#=20/100/200$, $d=8$, $\bsbb = (\{3\}^1, \{1.5\}^1, \{0\}^2, \{2\}^1,
\{0\}^3)$, $\Sigma_{ij}=\rho^{|i-j|}$ with $\rho=0.85$, $\sigma=2, 3, 5, 8$; the corresponding signal-to-noise variance ratio  is $8.21$, $3.65$, $1.31$, and $0.51$, respectively.\\
Before an algorithm is applied, the columns of a regression matrix are all normalized to have a squared $l_2$-norm
equal to the number of the observations; no centering is performed in these examples.

Each model is simulated 50 times, then, we measure the performance of each algorithm mainly by  test error and
 sparsity error.
The test error is characterized by the $40\%$ trimmed-mean of the scaled MSE (SMSE) on the test data,
where SMSE is $100\cdot (\sum_{i=1}^N ( \hat y_i - y_i )^2 /(N \sigma^2) - 1)$ defined for the test data.
(Medians of MSEs are mostly used~\cite{Tib, ZouHas}  to measure the performance from multiple runs,
but are not so stable for comparisons based on our experience.)
The sparsity error here is defined by the $40\%$ trimmed-mean of the following 50 percentages:
$100 \cdot |\{i:  \mbox{sgn}(\hat\beta_i)\neq \mbox{sgn}(\beta_i)\}|/d$, which represents the number of
inconsistent signs for each estimate compared to the true $\bsbb$.
We also summarized the proper zero percentages,
$100\% \cdot |\{i: \beta_i = 0, \hat\beta_i=0\}|/|\{i: \beta_i = 0\}|$, and the proper nonzero percentages,
$100\% \cdot |\{i: \beta_i \neq 0, \hat\beta_i\neq 0\}|/|\{i: \beta_i \neq 0\}|$ in the table as follows.
The numbers in parentheses are the standard errors of the trimmed means of SMSE, estimated by bootstrapping the  SMSE $500$ times as in~\cite{Zou}. The total computing time (in seconds) for each algorithm is also included.\\

\begin{table}[p]
\begin{center}
\small{
\begin{tabular}{|l|l||p{1.3cm}|p{1.4cm}|p{1.3cm}|p{1.3cm}||p{1.3cm}|p{1.3cm}|}

\hline
 \multicolumn{2}{|c||}{}& {Lasso} & {One-step SCAD} & {Hard-TISP} & {SCAD-TISP} & {eNet} & {Hybrid-TISP}\\

\hline
\multirow{4}{*}{EX1, $\sigma=2$}
    & {Test-err}  & \textbf{28.6}{\tiny{(3.6)}} & \textbf{25.3}{\tiny{(4.9)}} & \textbf{21.7}{\tiny{(3.7)}}& \textbf{18.2}{\tiny{(3.4)}} & \textbf{25.4}{\tiny{(3.5)}}& \textbf{15.9}{\tiny{(3.4)}}\\
    & {Spar-err} & \textbf{31.8} & \textbf{12.5}& \textbf{0} & \textbf{12.5} & \textbf{31.0} & \textbf{0} \\ \cline{2-8}
    & \emph{Prop-Z} & 50.8\% &91.2\%& 100\%& 89.5\% & 51.2\% & 100\% \\
    & \emph{Prop-NZ} & 100\% & 100\%& 100\%& 100\% & 100\% & 100\% \\

\hline
\multirow{4}{*}{EX1, $\sigma=3$}
    & {Test-err} & \textbf{27.8}{\tiny{(3.4)}} & \textbf{27.4}{\tiny{(3.5)}}& \textbf{25.9}{\tiny{(3.6)}}& \textbf{25.8}{\tiny{(3.5)}}& \textbf{23.4}{\tiny{(3.5)}}& \textbf{18.2}{\tiny{(4.0)}}\\
    & {Spar-err} & \textbf{30.7} & \textbf{16.7} & \textbf{5.5} & \textbf{12.5} & \textbf{31.5} & \textbf{3.8} \\ \cline{2-8}
    & \emph{Prop-Z} & 50.8\% & 80.0\% & 93.2\% & 92.0\% & 47.3\% & 94.1\% \\
    & \emph{Prop-NZ} &  100\% & 87.2\% & 100\% & 100\% & 100.0\% & 100.0\% \\

\hline
\multirow{4}{*}{EX1, $\sigma=5$}
    & {Test-err}  &  \textbf{23.0}{\tiny{(3.8)}}& \textbf{27.0}{\tiny{(2.5)}}& \textbf{22.3}{\tiny{(2.7)}}& \textbf{25.7}{\tiny{(3.1)}}& \textbf{18.4}{\tiny{(3.6)}}& \textbf{17.8}{\tiny{(3.5)}}\\
    & {Spar-err} & \textbf{32.0} & \textbf{25.0} & \textbf{12.5} & \textbf{25.0} & \textbf{31.5} & \textbf{12.5}   \\ \cline{2-8}
    & \emph{Prop-Z} & 50.4\% & 80.0\% & 91.6\% & 80.0\% & 48.6\% & 92.3\%  \\
    & \emph{Prop-NZ} &  86.2\% & 66.7\%& 85.1\% & 66.7\% & 100.0\% & 100\% \\

\hline
\multirow{4}{*}{EX1, $\sigma=8$}
    & {Test-err}  &  \textbf{15.4}{\tiny{(2.9)}}& \textbf{20.4 }{\tiny{(2.6)}}& \textbf{20.3}{\tiny{(2.9)}}& \textbf{17.1}{\tiny{(2.9)}}& \textbf{14.1}{\tiny{(2.8)}}& \textbf{11.3}{\tiny{(2.6)}}\\
    & {Spar-err} &  \textbf{31.3} & \textbf{30.5} & \textbf{25.0} & \textbf{32.8} & \textbf{37.5} & \textbf{30.6} \\ \cline{2-8}
    & \emph{Prop-Z} &  72.3\% & 80.0\% & 94.9\% & 80.0\% & 70.4\% & 94.7\% \\
    & \emph{Prop-NZ} &  66.7\% & 33.3\% & 49.5\% & 33.3\% & 66.7\% & 66.7\%\\

\hline
\hline

\multirow{4}{*}{EX2, $\sigma=2$}
    & {Test-err}  & \textbf{24.1}{\tiny{(2.1)}}& \textbf{28.7}{\tiny{(4.1)}} & \textbf{19.9}{\tiny{(2.8)}}& \textbf{20.8}{\tiny{(2.9)}}& \textbf{19.4}{\tiny{(3.5)}}& \textbf{14.3}{\tiny{(3.2)}} \\
    & {Spar-err} &  \textbf{31.0} & \textbf{32.0} & \textbf{12.5} & \textbf{12.5} & \textbf{30.8} & \textbf{12.5}   \\ \cline{2-8}
    & \emph{Prop-Z} &  60.0\% & 74.5\% & 80.0\% & 90.5\%  & 52.3\% & 80\% \\
    & \emph{Prop-NZ} & 100\% & 84.1\% & 100\% & 100\% & 100.0\% & 100.0\% \\

\hline
\multirow{4}{*}{EX2, $\sigma=3$}
    & {Test-err}  & \textbf{19.9}{\tiny{(3.3)}} & \textbf{29.1}{\tiny{(3.1)}} & \textbf{19.8}{\tiny{(2.2)}} & \textbf{20.2}{\tiny{(2.4)}} & \textbf{14.3}{\tiny{(2.9)}} & \textbf{14.1}{\tiny{(2.9)}}  \\
    & {Spar-err} & \textbf{30.7} & \textbf{29.8} & \textbf{16.1} & \textbf{16.3} & \textbf{28.6} & \textbf{17.3}\\ \cline{2-8}
    & \emph{Prop-Z} & 60.0\% & 80.0\% & 91.1\% & 91.9\% & 55.8\% & 80\%  \\
    & \emph{Prop-NZ} & 85.7\% & 66.7\% & 83.7\% & 66.7 \% & 100.0\% & 100.0 \% \\

\hline
\multirow{4}{*}{EX2, $\sigma=5$}
    & {Test-err}  & \textbf{13.9}{\tiny{(3.1)}} & \textbf{24.8}{\tiny{(2.2)}} & \textbf{16.6}{\tiny{(2.6)}} & \textbf{17.4}{\tiny{(2.7)}} & \textbf{9.7}{\tiny{(2.6)}} & \textbf{9.0}{\tiny{(2.7)}} \\
    & {Spar-err} &  \textbf{31.0} & \textbf{31.3} & \textbf{25.0}  & \textbf{25.0} & \textbf{30.0} & \textbf{25} \\ \cline{2-8}
    & \emph{Prop-Z} & 68.8\% & 80.0\% & 80.0\%  & 80.0\% & 52.9\% & 80\% \\
    & \emph{Prop-NZ} & 66.7\% & 48.2\% & 66.7\%  & 66.7\% & 100.0\% & 84.1\%\\

\hline
\multirow{4}{*}{EX2, $\sigma=8$}
    & {Test-err}  & \textbf{10.4}{\tiny{(2.6)}} & \textbf{16.3}{\tiny{(2.1)}} & \textbf{13.8}{\tiny{(2.8)}} & \textbf{15.6}{\tiny{(2.8)}} & \textbf{7.2}{\tiny{(2.2)}} & \textbf{6.9}{\tiny{(2.5)}} \\
    & {Spar-err} & \textbf{32.0} & \textbf{36.4} & \textbf{31.0} & \textbf{29.0} & \textbf{37.5} & \textbf{30.4}\\ \cline{2-8}
    & \emph{Prop-Z} & 71.0\% & 80.0\% & 92.1\% & 91.0\% & 49.3\% & 74.6\% \\
    & \emph{Prop-NZ} & 66.7\% & 33.3\% & 47.4\% & 33.3\% & 83.0\% & 66.7\%  \\

\hline
\hline
\multicolumn{2}{|c||}{Total computing time} &
284 & 599 & 596 & 610 & 2185 & 2234 \\
\hline
\end{tabular}\\
}


\end{center}
\caption{\small{Performance comparisons on the simulation data, in terms of
test error, sparsity error, proper sparsity, and proper nonsparsity -- all the numbers are 40\% trimmed-mean of the 50 simulations. Six methods are listed here: lasso (Soft-TISP), one-step SCAD, Hard-TISP, SCAD-TISP, elastic net (eNet), and Hybrid-TISP; the last two both have two regularization parameters. The last row gives the total computing time in seconds.}}
\label{table:simu1}
\end{table}

First, although Zou \& Li's one-step SCAD brings more sparsity than the lasso estimate (seen from the proper-sparsity and proper-nonsparsity), it is often the worst in terms of test error.
This is because  the one-step SCAD is indeed a weighted lasso method and the OLS estimate used for weight construction may not be trustworthy, if, say, there is large noise, or high correlation between some variables.
This phenomenon is  serious in Example 2 where the OLS estimate can be unstable and misleading.
Our Hard-TISP and SCAD-TISP clearly showed the remarkable \emph{parsimoniousness} brought by nonconvex penalties.
Instead of solving a $l_1$-constrained convex approximation as in the LLA method, our TISPs directly tackled the original nonconvex penalized regressions  and demonstrated better performance in both test-error and sparsity-error.
(In fact, we {doubt} if the $l_1$-based one-step SCAD is truly able to solve the  SCAD penalized regression, seen from the convex approximation in its derivation, and after comparing its estimate to the SCAD-TISP.) 
Hard-TISP and SCAD-TISP do not differ much here, which verifies the previous theoretical results regarding the hard-thresholding family in Section \ref{sectheo}.

Hard-TISP and SCAD-TISP achieve smaller test error than the lasso which may introduce extra bias when the signal-to-noise ratio is medium or high. Interestingly,  when the noise level is very high, the lasso (Soft-TISP) yields a more accurate estimate than the two.
This is in fact not so surprising. In predictive learning, to reduce the test error, when the noise is relatively large compared to the signal, it is also necessary to shrink the nonzero coefficients even if the true ones are far  from zero. In either hard- or SCAD-thresholding, there is basically no shrinkage offered for large nonzero coefficients, while the lasso does this by soft-thresholding (although the shrinkage amount is the same as the thresholding value).
Fortunately,  TISP still gives us good selection results and achieves parsimonious models. We can apply, for example, a second-time shrinkage  to the coefficients of the selected variables. Of course, a better strategy is to take into account these two concerns -- selection and shrinkage -- \emph{simultaneously} and \emph{adaptively} in building a model as probed in the next subsection.

\subsection{Hybrid-TISP for model selection and shrinkage}
\label{subsechybridtisp}
To deal with the low SNR problem, a promising approach is to modify the thresholding in Hard-TISP to include adaptive shrinkage for nonzero coefficients.
Motivated by the thresholding function of  ridge regression given by \eqref{ridgethfunc}, we propose a novel \emph{hybrid}-thresholding:
\begin{eqnarray}
\Theta(t;\lambda,\eta)=\begin{cases} 0, \quad \mbox{ if } |t|< \lambda\\ \frac{t}{1+\eta}, \quad \mbox{ if }  |t|\geq\lambda \end{cases}. \label{hybridthfunc}
\end{eqnarray}
The penalty constructed via  the mechanism introduced in Section \ref{subsecthpconst} is made up of two quadratic parts:
\begin{eqnarray}
P(\theta;\lambda, \eta)=\begin{cases} -\frac{1}{2} \theta^2 + \lambda |\theta|, \quad\mbox{ if } |\theta| < \frac{\lambda}{1+\eta}\\ \frac{1}{2} \eta \theta^2 +\frac{1}{2}\frac{\lambda^2}{1+\eta}, \quad\mbox{ if } |\theta| \geq \frac{\lambda}{1+\eta} \end{cases}. \label{hybridpenfunc}
\end{eqnarray}
We have seen the first quadratic part in the continuous hard-penalty \eqref{hardPcont} (which leads to the same solution as the discrete $l_0$-penalty); the second part resembles a ridge penalty. See Figure \ref{fighybrid} below. Note that the knots $\pm {\lambda}/{(1+\eta)}$ are  dependent on $\eta$, too.
Simple calculations show that  this $P$ satisfies the BCC (cf. \eqref{bcc}) with $\bsbH=\bsbI$, and Theorem \ref{conv}  holds. We can apply \eqref{tispvar} given an arbitrary design matrix.
The corresponding TISP (referred to as \textbf{Hybrid-TISP}) converges. The $\Theta$-equation \eqref{tispkkt} implies the nonzero components of a Hybrid-TISP estimate result from a partial ridge regression. This fact can be used in implementation when the maximum number of iterations allowed has been reached.

\begin{figure}[h!]
\begin{center}
\includegraphics[width=5in]{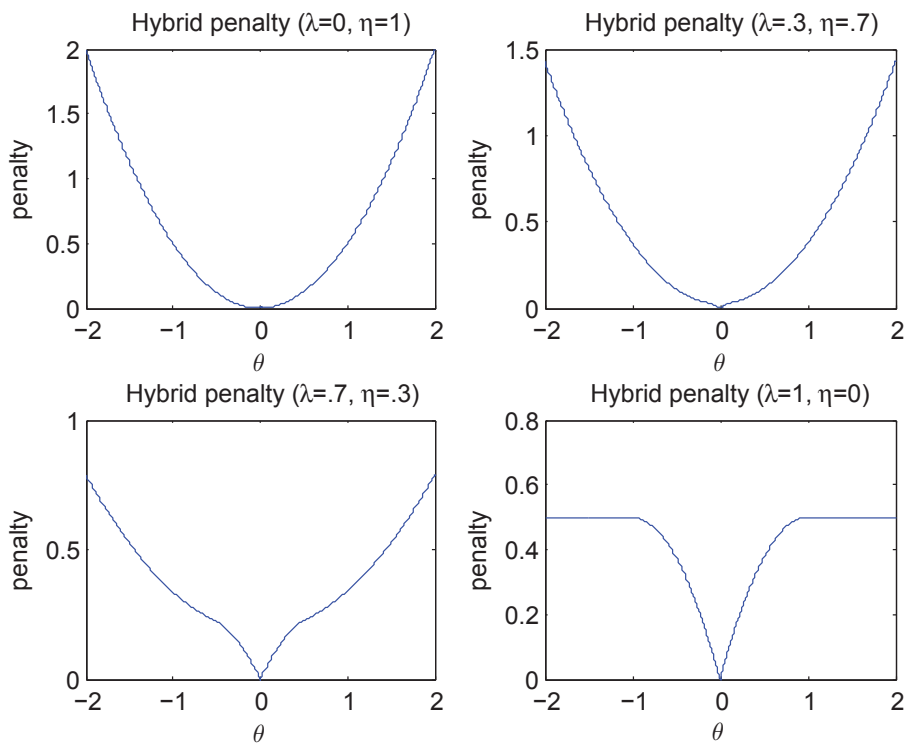}
\end{center}
\caption[Hybrid-penalty]{\small{The penalty defined by hybrid-thresholding. As $\lambda$ and $\eta$ vary, it takes the continuous hard-penalty and the ridge penalty as  extremes. }}
\label{fighybrid}
\end{figure}

Moreover, we have the following nonasymptotic result in parallel to Theorem \ref{selbyhard}. Recall that $k_0=\|\bsbX\|_2$, $\bsbSig^{(s)}=\bsbSig/n=\bsbX^T\bsbX/n$, $\mu = \mu_{\min}(\bsbSig_{nz,nz}^{(s)})$ and $\kappa \triangleq \underset{ i \in z}{\max} \| \bsbSig_{i,nz}^{(s)} \|_2 / \sqrt{ d_{nz}}$. Define $\iota\triangleq \min |(\bsbSig_{nz}+\eta\bsbI)^{-1}\bsbSig_{nz} \bsbb_{nz}|$, the minimum absolute value in the noiseless partial ridge estimate. Let $p_e$ be the probability of  Hybrid-TISP estimates having incorrect sparsity patterns, that is, for any $\hat \bsbb \in F$, there exists some $i$ or $j$ such that $\hat\bsbb_{z,i}\neq 0$ or $\hat\bsbb_{nz,j}=0$.
\begin{theorem}
\label{hybtispTH}
Assume $\mu>0$, and $\lambda,\eta$ are chosen such that $\kappa \leq \frac{\lambda}{\|\bsbb_{nz}\|_2\sqrt{ d_{nz}} } \frac{n\mu+\eta}{n\eta}$ and $\iota\geq \frac{\lambda}{k_0^2+\eta}$. Then
\begin{eqnarray}
p_e \leq 2d_z\varphi(M'')/M'' +2 d_{nz}\varphi(L'')/L'',\label{selbndhybridtisp}
\end{eqnarray}
where $M''=\frac{1}{\sqrt n \sigma}\left({\lambda} - \frac{ n \eta}{n\mu+\eta}\kappa  \| \bsbb_{nz}\|_2\sqrt {d_{nz}}\right)$, $L''
= \frac{n\mu+\eta}{\sqrt{n\mu}\sigma} \left(  \iota-\frac{\lambda}{k_0^2+\eta} \right)$.
\end{theorem}

Hybrid-TISP successfully offers both selection and shrinkage in estimating $\bsbb$. Before going into the numerical results, we summarize the traits of the design of Hybrid-TISP  as follows.
(a) Its penalty provides us a trade-off between the $l_0$-penalty and the $l_2$-penalty (ridge-penalty), and takes the two as extremes, from which we secure  selection and shrinkage simultaneously.
In particular, the selection is achieved by a penalty more like $l_0$ than $l_1$, seen from the penalty function, or the iterative thresholding.
(b) Hybrid-TISP avoids double shrinkage. Double shrinkage is a serious problem in the design of  naive elastic net~\cite{ZouHas} which simply adopts a  linear combination of the $l_1$-penalty and the $l_2$-penalty. However, the $l_1$-penalty also plays a role in shrinking the nonzero coefficients in addition to the $l_2$-penalty.
By contrast, Hybrid-TISP deals with the zeros and the nonzeros separately, by hard-thresholding and ridge-thresholding, respectively; there is no overlapping between them.
(c) We have two parameters, $\lambda$ and $\eta$, responsible for selection and shrinkage respectively.
One drawback of  the lasso is that it uses the same parameter to control both selection and shrinkage~\cite{Relaxo}.
Therefore, it may result in insufficient zeros even if the SNR is pretty high, as shown clearly in Table \ref{table:simu1}. Hybrid-TISP has $\lambda$, $\eta$  designed for  the two different purposes and can adapt to different sparsity and noise level.
(d) The TISP selecting and  shrinking  \emph{interplay} with each other  during the iteration till in the end we successfully achieve  selection/shrinkage balance in the final estimate. This is in contrast to the relaxed lasso~\cite{Relaxo} which treats selection and shrinkage as  separate steps in building a model.
(e) Finally, Hybrid-TISP is a very simple procedure to implement. It only involves multiplication and thresholding operations.

In the implementation of  Hybrid-TISP, an empirical parameter search is usually needed to determine the values of $\lambda$ and $\eta$, because running a grid search over the $(\lambda,\eta)$-space is a formidable task. We search along a couple of few one-dimensional solution paths including the $\lambda$-paths (with $\eta$ fixed) and the $\eta$-paths (with $\lambda$ fixed) to save computational cost. The optimal tuning parameter from the ridge regression path (corresponding to $\lambda=0$), denoted by $\eta^{(r)}$, is used as a reference  for $\eta$. Briefly, our search process generates and searches along some $\lambda$- and $\eta$-paths, compares the results from these searches, and then takes $(\lambda, \eta)$ to be the one minimizing the validation error. The concrete search  paths are as follows.
(i) $n>p$. Denote the OLS scale estimate by $\hat\sigma$. If $n/p<5$, or $n/p<10$ but $\hat\sigma>5$, we adopt the \emph{alternative} search strategy which has been shown to be fast and efficacious~\cite{Shethesis}:  fixing $\eta$ at $0.5\eta^{(r)}$, search along the $\lambda$-path to get an optimal solution (having the smallest validation error) at, say, $\lambda^{(o)}$; then  search along the $\eta$-path with $\lambda$  fixed at  $\lambda^{(o)}$. If $n/p>10$ and $\hat\sigma<5$, we only search over the $\lambda$-path with $\eta=0.05\eta^{(r)}$. In all remaining cases, we generate and search long two $\lambda$-paths with $\eta=0.5\eta^{(r)}$ and $0.05\eta^{(r)}$ respectively. (ii) $p>n$. We use the above alternative search starting with $0.5\eta^{(r)}$, and an additional search for  $\lambda$ with $\eta$ fixed at $0.05\eta^{(r)}$. Accordingly, 3 paths in total are generated in the large-$p$ situation. This simple empirical search does not cover the full parameter space but is more efficient than a grid search.
The results are reported in Table \ref{table:simu1}.
We also included the elastic net (eNet) in the experiments, which has two regularization parameters as well. Note that eNet generates and searches along 6 solution paths to tune the parameters~\cite{ZouHas}.

Seen from Table \ref{table:simu1}, Hybrid-TISP has amazing performance in both accuracy and sparsity.
We briefly summarize the story as follows.  When the noise level is low or medium, the value of $\lambda$ in the lasso is limited by the amount of shrinkage and thus gives insufficient sparsity. Large noise alleviates the problem but there is still much room for the improvement of test-error and sparsity-error because the amount of shrinkage may not equal to the thresholding value in the selection. The weighted lasso like the one-step SCAD has somewhat limited power  because the OLS estimate may be inaccurate and misleading for weight construction. Benefiting from the $l_2$-penalty, the eNet shows much better accuracy in the case of large noise and/or high correlation between the variables; nevertheless, the sparsity of the estimate may be  {seriously} hurt  when the ridge penalty must take control. And it seems possible to improve its test-error further by incorporating this sparsity in estimation. All of these problems can be resolved by Hybrid-TISP, which achieves the right balance between shrinkage and selection. Its test error is consistently lower than the eNet, and more importantly, Hybrid-TISP provides a parsimonious model as Hard-TISP.

\subsection{Large sample and large dimension experiments}
At the end of this section, we demonstrate the performance of TISP on large-$n$ data as well as on large-$p$ data. We modified the parameters in Example 1 and reran the simulations, where $\Sigma_{ij}=\rho^{|i-j|}$ with $\rho=.5$,  $\bsbb$ is appended with zeros given by $[3, 1.5, 0, 0, 2, 0, 0, \cdots, 0]^T$, $\sigma=2, 5$, and $n$, $d$ are not fixed anymore: in the large sample experiment, $d=8$, $n=40, 80, 200$ (corresponding to $5$ times, $10$ times, and $25$ times as large as $d$); in the large dimension experiment, $n=20$, $d=100, 200, 500$ (corresponding to $5$ times, $10$ times, and $25$ times as large as $n$). Table \ref{table:simu2} shows the simulation results of these different combinations of $n$ and $d$.
In both situations, the Hybrid-TISP path is preferable in terms of accuracy and sparsity.
Our conclusions are similar to the  findings summarized before.
Note that one-step SCAD uses the OLS estimate as the initial guess and thus is not included in the large-$p$ simulation. In fact, as an example of the adaptive lasso, it is most powerful in large samples with  small noise and low correlation between covariates, where the OLS estimate is accurate.
The elastic net is an improvement of the lasso and provides a good algorithm in predictive learning. However, despite having two regularization parameters, it does not improve much the sparsity of the lasso.
Hard- and SCAD-IPOD give significantly different solution paths than the above convex penalties. They may dramatically reduce the sparsity error and the test error, say, for large-$p$ sparse signals with moderate  noise. Both thresholdings fall into the hard-thresholding family which does not introduce much estimation bias for large coefficients. Interestingly, if our main concern is to reduce the test error in building a statistical model (which is the most frequently used tuning criterion in implementation), they are not always our best choices.   Indeed, it is more desirable to offer adaptive shrinkage to nonzero coefficient estimation to benefit from the bias-variance tradeoff. Hybrid-TISP is successful especially for the large-$p$ data because it does joint and adaptive selection and shrinkage.

\begin{table}[p]
\begin{center}
\small{
\begin{tabular}{|c|l||p{1.3cm}|p{1.4cm}|p{1.3cm}|p{1.3cm}||p{1.3cm}|p{1.3cm}|}

\hline
 \multicolumn{2}{|c||}{}& {Lasso} & {One-step SCAD} & {Hard-TISP} & {SCAD-TISP} & {eNet} & {Hybrid-TISP}\\


\hline
\multirow{2}{*}{$n=\textbf{40}$, $d=8$,}
    & {Test-err}  & \textbf{14.2}{\tiny{(2.8)}} & \textbf{10.6}{\tiny{(2.3)}} & \textbf{10.4}{\tiny{(2.7)}}& \textbf{9.2}{\tiny{(1.9)}} & \textbf{10.7}{\tiny{(2.7)}}& \textbf{8.2}{\tiny{(2.5)}}\\
    & {Spar-err} & \textbf{29.3} & \textbf{12.5}& \textbf{0.0} & \textbf{4.5} & \textbf{29.7} & \textbf{0.0} \\ \cline{2-8}
\multirow{2}{*}{ $\sigma=2$}
    & \emph{Prop-Z} & 53.0\% & 91.0 \%& 100\%& 100\% & 52.5\% & 100.0\% \\
    & \emph{Prop-NZ} & 100\% & 100\%& 100\%& 100\% & 100\% & 100\% \\
\hline
\multirow{2}{*}{$n=\textbf{40}$, $d=8$,}
    & {Test-err}  & \textbf{12.6}{\tiny{(3.0)}} & \textbf{15.2}{\tiny{(3.1)}} & \textbf{13.9}{\tiny{(2.4)}}& \textbf{13.3}{\tiny{(2.0)}} & \textbf{10.0}{\tiny{(2.9)}}& \textbf{12.3}{\tiny{(2.5)}}\\
    & {Spar-err} & \textbf{29.8} & \textbf{25.0}& \textbf{16.9} & \textbf{17.6} & \textbf{31.7} & \textbf{16.9} \\ \cline{2-8}
\multirow{2}{*}{$\sigma=5$}
    & \emph{Prop-Z} & 54.7\% & 66.4 \%& 94.1\%& 80\% & 52.0\% & 93.9\% \\
    & \emph{Prop-NZ} & 100\% & 83.0\%& 87.7\%& 83.3\% & 100.0\% & 100.0\% \\	
\hline
\multirow{2}{*}{$n=\textbf{80}$, $d=8$,}
    & {Test-err} & \textbf{9.2}{\tiny{(2.6)}} & \textbf{5.6}{\tiny{(2.0)}}& \textbf{5.7}{\tiny{(1.8)}}& \textbf{5.7}{\tiny{(1.9)}}& \textbf{8.7}{\tiny{(2.8)}}& \textbf{5.6}{\tiny{(1.8)}}\\
    & {Spar-err} & \textbf{44.3} & \textbf{3.5} & \textbf{0} & \textbf{4.7} & \textbf{42.9} & \textbf{0} \\ \cline{2-8}
\multirow{2}{*}{$\sigma=2$}
    & \emph{Prop-Z} & 29.2\% & 94.4\% & 100\% & 92.4\% & 31.3\% & 100\% \\
    & \emph{Prop-NZ} &  100\% & 100\% & 100\% & 100\% & 100.0\% & 100\% \\
\hline
\multirow{2}{*}{$n=\textbf{80}$, $d=8$,}
    & {Test-err} & \textbf{8.8}{\tiny{(2.5)}} & \textbf{10.0}{\tiny{(2.3)}}& \textbf{5.8}{\tiny{(1.8)}}& \textbf{6.8}{\tiny{(1.9)}}& \textbf{7.6}{\tiny{(2.2)}}& \textbf{5.2}{\tiny{(1.7)}}\\
    & {Spar-err} & \textbf{44.9} & \textbf{17.8} & \textbf{4.6} & \textbf{17.9} & \textbf{37.5} & \textbf{3.8} \\ \cline{2-8}
\multirow{2}{*}{$\sigma=5$}
    & \emph{Prop-Z} & 29.2\% & 80\% & 94.9\% & 80\% & 40.0\% & 95.1\% \\
    & \emph{Prop-NZ} &  100\% & 100\% & 100\% & 100\% & 100\% & 100.0\% \\
\hline
\multirow{2}{*}{$n=\textbf{200}$, $d=8$,}
    & {Test-err} & \textbf{3.2}{\tiny{(1.6)}} & \textbf{1.7}{\tiny{(1.4)}}& \textbf{1.8}{\tiny{(1.4)}}& \textbf{1.8}{\tiny{(1.4)}}& \textbf{2.3}{\tiny{(1.5)}}& \textbf{1.8}{\tiny{(1.4)}}\\
    & {Spar-err} & \textbf{43.4} & \textbf{3.2} & \textbf{3.8} & \textbf{2.9} & \textbf{29.8} & \textbf{3.8} \\ \cline{2-8}
\multirow{2}{*}{$\sigma=2$}
    & \emph{Prop-Z} & 30.5\% & 94.8\% & 93.9\% & 95.3\% & 52.4\% & 93.9\% \\
    & \emph{Prop-NZ} &  100\% & 100\% & 100\% & 100\% & 100\% & 100\% \\

\hline
\multirow{2}{*}{$n=\textbf{200}$, $d=8$,}
    & {Test-err} & \textbf{3.2}{\tiny{(1.6)}} & \textbf{2.5}{\tiny{(1.5)}}& \textbf{2.2}{\tiny{(1.4)}}& \textbf{2.7}{\tiny{(1.6)}}& \textbf{3.0}{\tiny{(1.4)}}& \textbf{1.9}{\tiny{(1.3)}}\\
    & {Spar-err} & \textbf{37.5} & \textbf{17.8} & \textbf{12.5} & \textbf{3.3} & \textbf{29.8} & \textbf{12.5} \\ \cline{2-8}
\multirow{2}{*}{$\sigma=5$}
    & \emph{Prop-Z} & 40.0\% & 71.1\% & 92.4\% & 98.0\% & 52.3\% & 80.0\% \\
    & \emph{Prop-NZ} &  100\% & 100\% & 100\% & 100\% & 100\% & 100\% \\
\hline
\hline
\multirow{2}{*}{$n=20$, $d=\textbf{100}$,}
    & {Test-err}  &  \textbf{83.2}{\tiny{(7.9)}}& ------ & \textbf{68.3}{\tiny{(10.6)}}& \textbf{45.1}{\tiny{(8.8)}}& \textbf{82.9}{\tiny{(7.4)}}& \textbf{64.3}{\tiny{(10.2)}}\\
    & {Spar-err} & \textbf{7.9} & ------ & \textbf{1.3} & \textbf{1.0} & \textbf{7.5} & \textbf{1.3}   \\ \cline{2-8}
\multirow{2}{*}{$\sigma=2$}
    & \emph{Prop-Z} & 100\% & ------ & 100\% & 100.0\% & 100.0\% & 98.5\%  \\
    & \emph{Prop-NZ} &  92.0\% & ------ & 99.3 \% & 99.6 \% & 92.2\% & 100.0\% \\
\hline
\multirow{2}{*}{$n=20$, $d=\textbf{100}$,}
    & {Test-err}  &  \textbf{50.7}{\tiny{(3.7)}}& ------ & \textbf{53.0}{\tiny{(7.1)}}& \textbf{50.0}{\tiny{(6.8)}}& \textbf{51.0}{\tiny{(4.7)}}& \textbf{39.9}{\tiny{(6.1)}}\\
    & {Spar-err} & \textbf{7.5} & ------ & \textbf{3.5} & \textbf{2.3} & \textbf{7.6} & \textbf{2.4}   \\ \cline{2-8}
\multirow{2}{*}{$\sigma=5$}
    & \emph{Prop-Z} & 93.9\% & ------ & 99.6\% & 99.7\% & 93.2\% & 98.7\%  \\
    & \emph{Prop-NZ} &  52.8\% & ------ & 49.6\% & 33.3\% & 66.7\% & 66.7\% \\
\hline
\multirow{2}{*}{$n=20$, $d=\textbf{200}$,}
    & {Test-err}  &  \textbf{121.2}{\tiny{(15.0)}}& ------ & \textbf{98.0}{\tiny{(12.0)}}& \textbf{92.3}{\tiny{(16.3)}}& \textbf{121.6}{\tiny{(13.9)}}& \textbf{103.9}{\tiny{(12.6)}}\\
    & {Spar-err} &  \textbf{4.7} & ------ & \textbf{1.0} & \textbf{0.8} & \textbf{4.2} & \textbf{0.7} \\ \cline{2-8}
\multirow{2}{*}{$\sigma=2$}
    & \emph{Prop-Z} &  95.4\% & ------ & 99.9\% & 99.9\% & 95.7\% & 99.9\% \\
    & \emph{Prop-NZ} &  100\% & ------ & 84.6\% & 85.0\% & 100.0\% & 66.7\%\\
\hline

\multirow{2}{*}{$n=20$, $d=\textbf{200}$,}
    & {Test-err}  &  \textbf{57.2}{\tiny{(5.8)}}& ------ & \textbf{58.6}{\tiny{(6.8)}}& \textbf{59.8}{\tiny{(7.5)}}& \textbf{54.1}{\tiny{(5.8)}}& \textbf{47.7}{\tiny{(6.5)}}\\
    & {Spar-err} &  \textbf{3.5} & ------ & \textbf{1.5} & \textbf{1.8} & \textbf{4.6} & \textbf{1.3} \\ \cline{2-8}
\multirow{2}{*}{$\sigma=5$}
    & \emph{Prop-Z} &  97.5\% & ------ & 99.2\% & 99.3\% & 95.8\% & 99.0\% \\
    & \emph{Prop-NZ} &  46.0\% & ------ & 53.3\% & 51.6\% & 66.7\% & 66.7\%\\
\hline
\multirow{2}{*}{$n=20$, $d=\textbf{500}$,}
    & {Test-err}  &  \textbf{173.5}{\tiny{(19.1)}}& ------ & \textbf{117.9}{\tiny{(14.6)}}& \textbf{118.1}{\tiny{(16.3)}}& \textbf{181.1}{\tiny{(20.1)}}& \textbf{118.2}{\tiny{(14.5)}}\\
    & {Spar-err} & \textbf{2.5} & ------ & \textbf{0.4} & \textbf{0.3} & \textbf{2.1} & \textbf{0.4}   \\ \cline{2-8}
\multirow{2}{*}{$\sigma=2$}
    & \emph{Prop-Z} & 97.6\% & ------ & 100.0\% & 100\% & 97.9\% & 99.9\%  \\
    & \emph{Prop-NZ} &  100\% & ------ & 66.7\% & 66.7\% & 88.4\% & 66.7\% \\
\hline
\multirow{2}{*}{$n=20$, $d=\textbf{500}$,}
    & {Test-err}  &  \textbf{72.9}{\tiny{(7.5)}}& ------ & \textbf{66.6}{\tiny{(6.5)}}& \textbf{67.5}{\tiny{(7.0)}}& \textbf{67.3}{\tiny{(7.2)}}& \textbf{60.1}{\tiny{(6.5)}}\\
    & {Spar-err} & \textbf{1.2} & ------ & \textbf{0.7} & \textbf{0.7} & \textbf{1.4} & \textbf{1.1}   \\ \cline{2-8}
\multirow{2}{*}{$\sigma=5$}
    & \emph{Prop-Z} & 99.2\% & ------ & 99.5\% & 99.6\% & 99.0\% & 99.3\%  \\
    & \emph{Prop-NZ} &  33.3\% & ------ & 50.6\% & 50.0\% & 33.3\% & 51.3\% \\
\hline
\end{tabular}\\
}


\end{center}
\caption{\small{Performance comparisons on the simulation data with large sample size or dimensionality, in terms of test error, sparsity error, proper sparsity, and proper nonsparsity over 50 simulations. Six methods are listed here: lasso (Soft-TISP), one-step SCAD, Hard-TISP, SCAD-TISP, elastic net (eNet), and Hybrid-TISP; the last two both have two regularization parameters.}}
\label{table:simu2}
\end{table}

\section{Real Data}
Hybrid-TISP was applied to a real  prostate dataset which was used by Tibshirani~\cite{Tib}.
The prostate data have $97$ observations and $9$ clinical measures.
In this example, unlike \cite{Tib}, we take the log(cancer volume) ({\tt{lcavol}}) as the response variable and consider a full quadratic model; the $43$ predictors are $8$ main effects, $7$ squares, and $28$ interactions of  eight original variables --- {\tt lweight}, {\tt age}, {\tt lbph}, {\tt svi}, {\tt lcp}, {\tt gleason},  {\tt pgg45}, and {\tt lpsa}, where {\tt svi} is binary.
The lasso does not give stable and accurate results for this example due to the existence of many highly correlated predictors.

\begin{figure}[h]
\begin{center}
\includegraphics[width=5in]{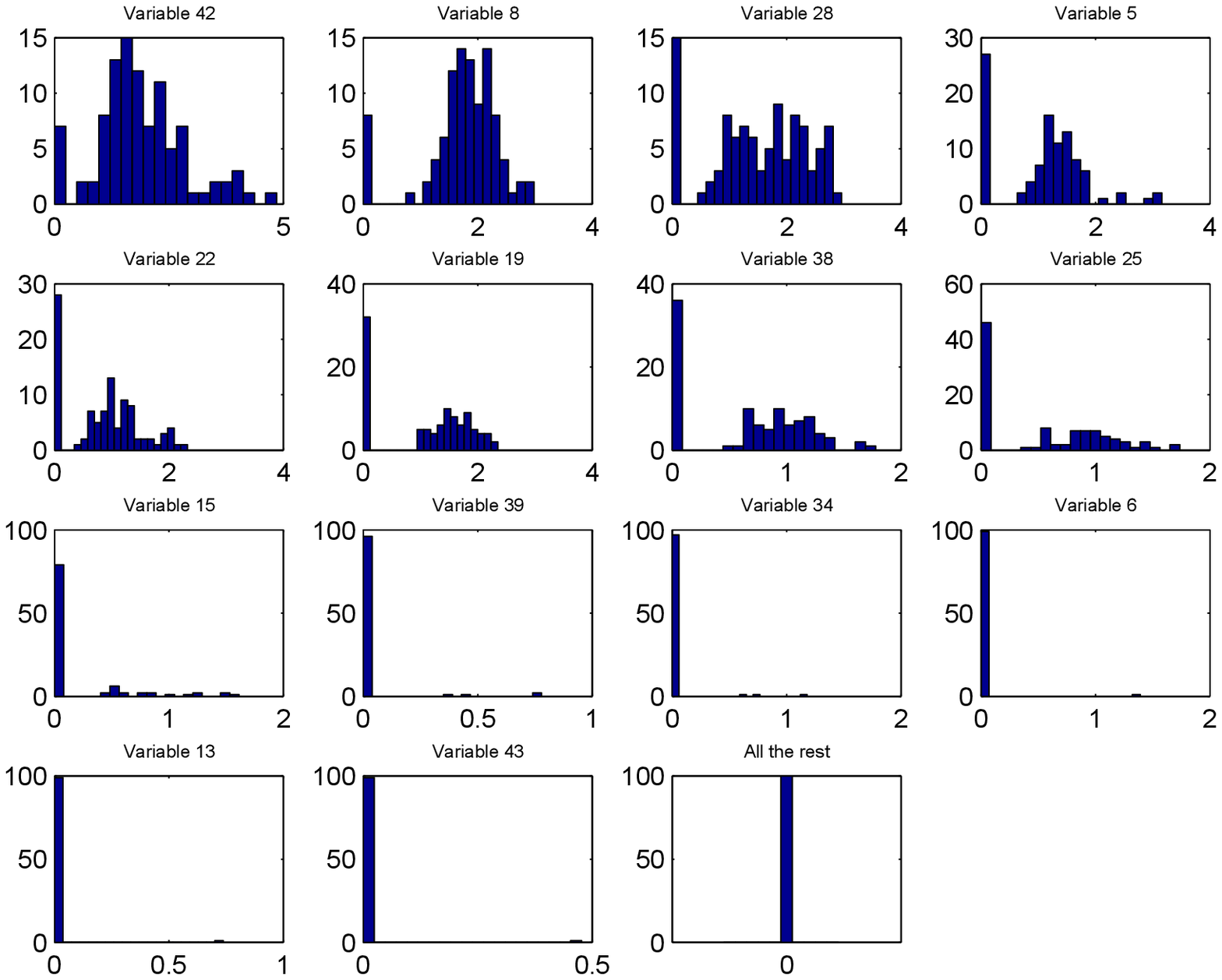}
\end{center}
\caption[Hybrid-TISP coefficient histograms]{\small{Histograms of the $100$ Hybrid-TISP coefficient estimates for all of the $43$ predictors in the prostate example.}}
\label{fighists}
\end{figure}

\begin{figure}[h]
\begin{center}
\includegraphics[width=4in]{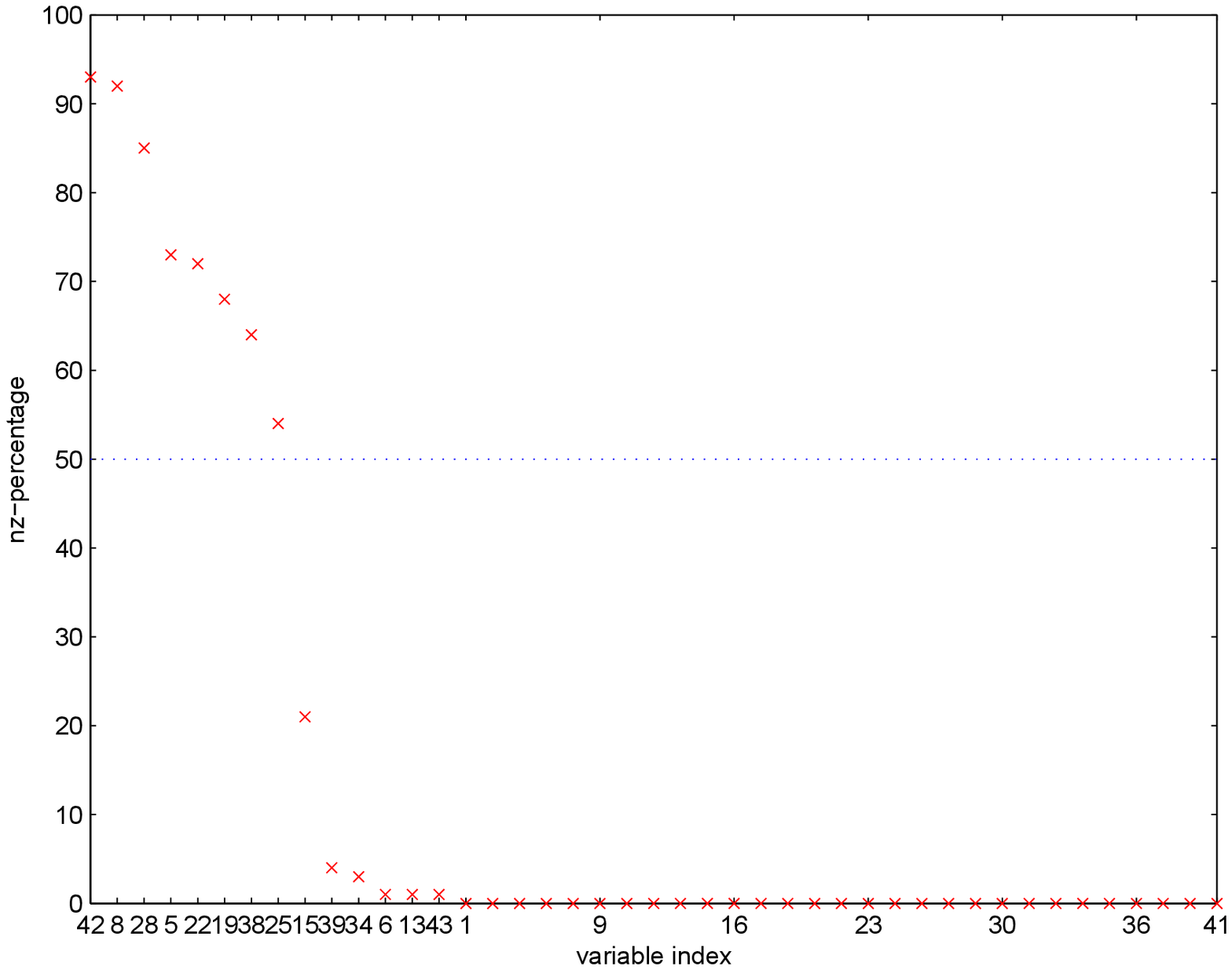}
\end{center}
\caption[Non-zero Hybrid-TISP coefficient proportions]{\small{Proportions of the Hybrid-TISP coefficients being nonzero over the $100$ bootstrap replications in the prostate example. Eight variables have nonzero coefficient estimates more frequently  than zero estimates.
}}
\label{figprops}
\end{figure}

The regularization parameters of Hybrid-TISP were tuned by leave-one-out cross-validation.
To identify the relevant variables in a trustworthy way, nonparametric bootstrap resampling was used with $B=100$.
For every bootstrap dataset, after standardizing the predictors, we apply  Hybrid-TISP with fixed regularization parameters  tuned for the original dataset.
Figure \ref{figprops} shows the percentages of the bootstrap coefficient estimates being nonzero  over the $100$ replications for all  the $43$ predictors.
The histograms are plotted in Figure \ref{fighists}. It is easy to see that $8$ variables are much more significant than the others. In fact, these are exactly the variables selected by Hybrid-TISP on the original data.
A more careful examination shows that they  appear (jointly) $36$ times in the selected models, the top visited octuple in  bootstrapping.
These variables fall into two groups with similar patterns: (I) $\{x_5, x_{19}, x_{25}, x_{38}\}$, i.e.,  \{{\tt lcp},   {\tt lweight*lcp},   {\tt age*lcp},   {\tt gleason*lcp}\}; (II) $\{x_8, x_{22}, x_{28}, x_{42}\}$, i.e., \{{\tt lpsa}, {\tt lweght*lpsa}, {\tt age*lpsa}, {\tt gleason*lpsa}\}. The within-group correlations are very high, $>.98$ for  Group (I), and $>.93$ for Group (II).
Furthermore, an interesting feature is that for any of the eight  variables selected by Hybrid-TISP, the other  three in the same group are  most correlated  with it among  $42$ predictors.

\section{Discussion}
We have proposed the thresholding-based iterative selection procedures for solving nonconvex penalized regressions.
In fact, people have long before noticed the weakness of the convex $l_1$-constraint (or the soft-thresholding) in wavelets and have designed many different forms of nonconvex penalties to increase model sparsity and accuracy. But for a nonorthogonal regression matrix, there is great difficulty  in both investigating the  performance in theory and solving the problem  in computation. TISP provides a simple and efficient way to tackle this.

Somewhat different than  other studies, we started from  thresholding rules rather than penalty functions.
Indeed, there is  a universal connection between them.
But a drawback of the latter is its non-unique form: different penalties may result in the same estimator and the same thresholding. The main contribution of this paper is the study of a class of $\Theta$-estimators satisfying \eqref{tispkkt}, which can be naturally computed by TISP, and are associated with penalized regressions. With a carefully designed thresholding rule, we obtained a good estimator for model selection and shrinkage.
Starting from $\Theta$ greatly facilitated the computation and the analysis.
In fact, some penalty designs may even have a better explanation from $\Theta$, or equivalently, the $\psi$-function --- for example, the SCAD-penalty (recall that it is defined by its derivative) seems to originate from Hampel's three-part redescending  $\psi$. Conversely, we can use TISP to compute $M$-estimators in robust statistics as described in Section \ref{subsecrelwork}.

Using a thresholding rule in the hard-thresholding family, TISP gives good selection results. Our novel Hybrid-TISP, accomplishing a fusion between $l_0$-penalty and $l_2$-penalty based on the hard-thresholding and the ridge thresholding, shows superior performance and beats the commonly used methods in both test-error and sparsity. The hybrid penalty function \eqref{hybridpenfunc}
 may look a bit odd, but is quite natural and simple from the point of view of thresholding; see \eqref{hybridthfunc}. It is worth mentioning that in contrast to~\cite{Gao, AntFan, antrev}, where more than one tuning parameter is considered a drawback and unnecessity, we believe a good procedure should have two explicit regularization parameters to control and balance selection and shrinkage.

We assume the penalty function $P$ is dependent on $\bsbb$ and $\lambda$ only. Therefore the iterative weighting, substituting the nonnegative garrote~\cite{Gao} for $\Theta$ in TISP, is not covered by the studies in this paper. In fact, with $\bsbb$ involved in $P$, it might be difficult to optimize in the second step of the mechanism introduced in Section \ref{secmotiv}.

The solution path associated with a nonconvex penalty is generally not continuous in $\lambda$. For example, even for the transformed $l_1$-penalty in Example 3 which is differentiable to any order on $(0,+\infty)$, the solution path still has no $\lambda$-continuity practically.
Hence a pathwise algorithm is not appropriate here. Empirically, using a zero estimate as the start in nonconvex TISPs works pretty well. We conjecture that it leads to an estimate with some least norm property.
Take  Hard-IPOD as an example: this roughly means that we were looking for the local minimum of the $l_0$-penalized regression that is closest to zero in building a parsimonious model.

The generalization of TISP to GLM seems straightforward; we will investigate this topic in the next paper.
TISP fits perfectly into the Accelerated Annealing~\cite{Shethesis} and thus can be used in the generic sparse regression with customizable sparsity patterns, such as the supervised clustering problem.
Other future studies include developing some acceleration techniques for TISP  (like the relaxation and asynchronous updating~\cite{Shethesis}) and deriving some risk oracles in theory.

\section*{Acknowledgements}
The author is grateful to the two anonymous referees, and especially the associate editor, for careful comments and  useful suggestions.
Most of this paper is based on a previous technical report~\cite{SheTisp}, supported by NSF grant DMS-0604939. The author would like to thank Art Owen for his valuable guidance.

\appendix
\section{Proofs}

\subsection{Proofs of Theorem \ref{conv}, Proposition \ref{convstat}, and Proposition \ref{optimality}}
\label{appproof1}

%
%

Let's consider the orthogonal case first.
Define $Q(\bsbg)= \| \bsbg -\bsba\|_2^2/2 + P(\bsbg;\lambda)$, where $\bsba$ is a known vector. Let $\bsbg_o=\arg \min Q(\bsbg)$. By the construction of $P$ and Proposition \ref{uniqsol}, $\bsbg_o$ satisfies $\bsbg_o-\bsba + s(\bsbg_o;\lambda)=\bsb{0}$.

\begin{eqnarray*}
Q(\bsbg_o+\bsb{h})-Q(\bsbg_o)&=&\frac{1}{2}\|\bsbg_o+\bsb{h}-\bsba\|_2^2-\frac{1}{2}\|\bsbg_o-\bsba\|_2^2
+P(\bsbg_o+\bsb{h}; \lambda)-P(\bsbg_o;\lambda)\\
&=&\frac{1}{2}\|\bsb{h}\|_2^2+<\bsb{h},\bsbg_o-\bsba>+P(\bsbg_o+\bsb{h};\lambda)-P(\bsbg_o;\lambda)\\
&=&\frac{1}{2}\|\bsb{h}\|_2^2+\left( P(\bsbg_o+\bsb{h};\lambda)-P(\bsbg_o;\lambda) -<\bsb{h},\bsb{s}> \right)\\
&\geq&\frac{1}{2}\|\bsb{h}\|_2^2 - \frac{1}{2} \bsb{h}^T \bsbH \bsb{h} = \frac{1}{2}\bsb{h}^T(\bsbI-\bsbH)\bsb{h}.
\end{eqnarray*}
This inequality is due to the BCC \eqref{bcc}. On the other hand, we know
\begin{eqnarray*}
Q(\bsbg_o+\bsb{h})-Q(\bsbg_o) \geq 0.
\end{eqnarray*}
In summary, we get
\begin{eqnarray}
Q(\bsbg_o+\bsb{h})-Q(\bsbg_o) \geq \frac{1}{2} \bsb{h}^T\bsbA\bsb{h},  \label{orthosolbnd}
\end{eqnarray}
for both  $\bsbA=\bsbI-\bsbH$ and $\bsbA=\bsb{0}$; formally, we write $\bsbA=(\bsbI-\bsbH)\vee\bsb{0}$. Note that \eqref{orthosolbnd} is a global result  for \emph{any} $\bsb{h}$.

Now look at the TISP. Recall the $g$ in \eqref{gdef} is
\begin{eqnarray*}
g(\bsbb,\bsbg)=\frac{1}{2}\|\bsbX\bsbg-\bsby\|_2^2 + P(\bsbg;\lambda)+\frac{1}{2} (\bsbg-\bsbb)^T(\bsbI-\bsbSig) (\bsbg-\bsbb).
\end{eqnarray*}
Then given $\bsbb$, we can write $g$ as
$$
g(\bsbb, \bsbg) = \frac{1}{2} \| \bsbg-((\bsbI-\bsbSig)\bsbb+\bsbX^T\bsby)\|_2^2+C(\bsbX,\bsby,\bsbb),
$$
and apply \eqref{orthosolbnd} with $\bsba=(\bsbI-\bsbSig)\bsbb+\bsbX^T\bsby$,
\begin{eqnarray}
g(\bsbb,\bsbg_o(\bsbb)+\bsb{h})-g(\bsbb,\bsbg_o(\bsbb)) \geq \frac{1}{2} \bsbh^T ((\bsbI-\bsbH)\vee \bsb{0}) \bsbh, \quad \forall \bsbh. \label{tispoptineq}
\end{eqnarray}
Correspondingly, for the TISP iterates $\bsbb^{(j)}$, we have
\begin{eqnarray*}
&&f(\bsbb^{(j+1)})+\frac{1}{2}(\bsbb^{(j+1)}-\bsbb^{(j)})^T(\bsbI-\bsbSig)(\bsbb^{(j+1)}-\bsbb^{(j)}) =  g(\bsbb^{(j)},\bsbb^{(j+1)}) \\
&\leq& g(\bsbb^{(j)},\bsbb^{(j)})-\frac{1}{2}(\bsbb^{(j+1)}-\bsbb^{(j)})^T ((\bsbI-\bsbH)\vee \bsb{0})(\bsbb^{(j+1)}-\bsbb^{(j)})\\
&=&f(\bsbb^{(j)})-\frac{1}{2}(\bsbb^{(j+1)}-\bsbb^{(j)})^T ((\bsbI-\bsbH)\vee \bsb{0})(\bsbb^{(j+1)}-\bsbb^{(j)}).
\end{eqnarray*}
That is,
\begin{equation}
f(\bsbb^{(j)})- f(\bsbb^{(j+1)})\geq \frac{1}{2}(\bsbb^{(j+1)}-\bsbb^{(j)})^T ((\bsbI-\bsbH)\vee \bsb{0}+\bsbI-\bsbSig)(\bsbb^{(j+1)}-\bsbb^{(j)}). \label{appasympreg}
\end{equation}
Now \eqref{optfval} and \eqref{asympreg} can be obtained after simple calculations.\\

As for Proposition \ref{convstat}, let $\bsbb^{(j_k)}\rightarrow\bsbb^*$ as $k\rightarrow\infty$. Under the condition $\mu_{\max}(\bsbSig)< 1\vee (2-\mu_{\max}(\bsbH))$, Theorem \ref{conv} states that
\begin{eqnarray*}
\| \bsbb^{(j_k+1)}-\bsbb^{(j_k)}\|_2^2 \leq (f(\bsbb^{(j_k)})- f(\bsbb^{(j_k+1)}))/C \leq (f(\bsbb^{(j_k)})- f(\bsbb^{(j_{k+1})}))/C\rightarrow 0.
\end{eqnarray*}
That is, $\Theta((\bsbI-\bsbSig)\bsbb^{(j_k)}+\bsbX^T \bsby;\lambda)-\bsbb^{(j_k)}\rightarrow 0$. Therefore, $\bsbb^*$ is a fixed point of TISP.\\

Finally, we prove Proposition \ref{optimality}. Noticing that $\bsbg_o(\bsbb^*)=\bsbb^*$, we get the following inequality from \eqref{tispoptineq}
\begin{eqnarray*}
g(\bsbb^*,\bsbb^*+\bsbh)-g(\bsbb^*,\bsbb^*)\geq \frac{1}{2} \bsbh^T ((\bsbI-\bsbH) \vee \bsb{0}) \bsbh, \quad \forall \bsbh.
\end{eqnarray*}
Since $g(\bsbb^*,\bsbb^*)=f(\bsbb^*)$,
\begin{eqnarray*}
&&f(\bsbb^*+\bsbh)+\frac{1}{2} \bsbh^T (\bsbI-\bsbSig) \bsbh \geq f(\bsbb^*)+ \frac{1}{2} \bsbh^T ((\bsbI-\bsbH) \vee \bsb{0}) \bsbh, \quad \forall \bsbh \\
&\Rightarrow&
f(\bsbb^*+\bsbh)-f(\bsbb^*) \geq  \frac{1}{2} \bsbh^T ((\bsbI-\bsbH) \vee \bsb{0}+\bsbSig-\bsbI) \bsbh, \quad \forall \bsbh.
\end{eqnarray*}
Therefore, if $\mu(\bsbSig)\geq\mu_{\max}(\bsbH)$, $\bsbb^*$ is a global minimizer of $f$.

\subsection{Proofs  of Theorem \ref{selbytisp} and Theorem \ref{selbyhard}}
\label{appproof2}
These theorems have all been essentially proved in~\cite{Shethesis}. We provide a self-contained proof as follows.
All inequalities and the absolute value `$||$' are understood
in the componentwise sense.
Assume, for the moment, $\bsbX$ has been column-normalized such that the diagonal entries of $\bsbSig=\bsbX^T\bsbX$ are all 1.
Let $\bsbSig_I=\bsbX_I^T\bsbX_I$, $\bsbSig_{I, I'}=\bsbX_I^T\bsbX_{I'}$ for any index sets $I, I'$.
\begin{lemma}
\label{dawl:lemma1}
Assume $\bsbSig_{nz}$ is nonsingular.
The TISP estimate $\hat\bsbb$ satisfies the following equations
\begin{equation}
\bsbS_z\hat \bsbb_z =   (\bsbX_z^{T} - \bsbSig_{z, nz} \bsbSig_{nz}^{-1} \bsbX_{nz}^T ) \bsbeps +
\tau \bsbSig_{z, nz} \bsbSig_{nz}^{-1} \widetilde{{\mbox{sgn}}}(\hat\bsbb_{nz}; {\lambda}/{k_0^2}) - \tau
\widetilde{{\mbox{sgn}}}(\hat\bsbb_z; {\lambda}/{k_0^2})   \label{dawl:(A2.8)}
\end{equation}
\begin{equation}
\hat \bsbb_{nz} = \bsbb_{nz} + \bsbSig_{nz}^{-1} ( \bsbX_{nz}^{T} \bsbeps - \tau \widetilde{{\mbox{sgn}}}(\hat\bsbb_{nz}; {\lambda}/{k_0^2})
) -\bsbSig_{nz}^{-1} \bsbSig_{z, nz}^T \hat\bsbb_z\label{dawl:(A2.9)}
\end{equation}
where  $\bsbS_z=\bsbSig_{z} - \bsbSig_{z, nz} \bsbSig_{nz}^{-1} \bsbSig_{nz, z}$.
\end{lemma}
This can be obtained directly from \eqref{kktgsgn}.

\begin{lemma}
\label{dawl:lemma2}
Let $\bsb{z} \sim N(\bsb{0}, \bsb{D}_{d\times d})$, $\bsb{z}'
\sim N(\bsb{0}, \bsb{\Lambda}_{d\times d})$, where $\bsb{D}$ is a diagonal matrix. Assume the diagonal
entries of $\bsb{\Lambda}$, denoted by $\mbox{diag}(\bsb{\Lambda})$, are the same as those of $\bsb{D}$,
i.e., $\mbox{diag}(\bsb{\Lambda})=\mbox{diag}(\bsb{D})$. Then
\begin{eqnarray}
P(\max |\bsb{z}'| \geq c) \leq P(\max |\bsb{z}| \geq c),   \label{dawl:(A2.10)}
\end{eqnarray} for any $c$.
\end{lemma}
This is clear from \v{S}id\'{a}k's classical result~\cite{Sid} in 1967.\\

To prove Theorem \ref{selbytisp}, let ${\bsbX_{z}'}^{T} =\bsbX_{z}^{T} - \bsbSig_{z, nz} \bsbSig_{nz}^{-1} \bsbX_{nz}^T$ and define
\begin{eqnarray*}
A&\triangleq & \left\{\left| \bsbX_{z} '
\bsbeps + \tau \bsbSig_{z, nz} \bsbSig_{nz}^{-1} \widetilde{{\mbox{sgn}}}(\hat\bsbb_{nz}) \right| \leq \tau \right\}\\
V&\triangleq &\left\{ \left|\bsbSig_{nz}^{-1} \bsbX_{nz}^{T} \bsbeps \right| +  \tau \left| \bsbSig_{nz}^{-1}
\bsb{s} \right| < \left| \bsbb_{nz} \right| {\mbox{ for any }} \bsb{s} {\mbox{ satisfying }}
|\bsb{s} |\leq 1 \right\}.
\end{eqnarray*}
Clearly $1-p_s\leq P(A^c\cup V^c)\leq P(A^c)+P(V^c)$.

Let $\bsbSig_{z, nz}=[\bsb{v}_{1},\cdots\bsb{v}_{d_{z}}]^T$, then $\|\bsb{v}_i\|_{2} \leq
\kappa \sqrt{d_{nz}}$. So $|\bsb{v}_i^T \bsbSig_{nz}^{-1} \widetilde{{\mbox{sgn}}}(\hat\bsbb_{nz}) | \leq
\|\bsb{v}_i\|_{2} \cdot \|\bsbSig_{nz}^{-1}\|_{2} \cdot \| \widetilde{{\mbox{sgn}}}(\hat\bsbb_{nz})\|_{2} \leq \kappa
d_{nz} / \mu$ and $\left| \bsbSig_{z, nz} \bsbSig_{nz}^{-1} \widetilde{{\mbox{sgn}}}(\hat\bsbb_{nz}) \right| \leq
\kappa d_{nz} / \mu$.
It follows that
$ P(A^c)\leq P\left(\left\{\max\left| {\bsbX_{z} '}^T \bsbeps \right| \geq (1-\kappa d_{nz} / \mu)\tau\right\} \right).$

Define $\bsbeps_{1}'={\bsbX_{z} '}^T \bsbeps \in \Re^{d_{z}}.$ Note that ${\bsbX_{z} '}^T {\bsbX_{z} '} = \bsbS_{z}$. Thus $\bsbeps_{1}'\sim N(\bsb{0}, \sigma^2\bsbS_{z})$. Since $\mbox{diag}(\bsbSig_{z})=\bsb{1}$ and
$\mbox{diag}(\bsbSig_{z, nz} \bsbSig_{nz}^{-1} \bsbSig_{z, nz}^T)=\left[ \bsb{v}_i^T \bsbSig_{nz}^{-1} \bsb{v}_i
\right] \geq \bsb{0}$, $\mbox{diag}(\bsbS_{z})\leq\bsb{1}$. Lemma \ref{dawl:lemma2} states that
\begin{eqnarray*}
P(A^c) &\leq& P\left(\max |\bsbeps_{1}''| \cdot \sigma \geq \tau(1-\frac{\kappa
d_{nz}}{\mu})\right) \\
&\leq& d_{z} P\left(|\epsilon_{1, i}''| \geq \tau (1-\frac{\kappa d_{nz}}{\mu})
\frac{1}{\sigma}\right)\equiv 2d_{z} \Phi ([M, +\infty)),
\end{eqnarray*} where $\bsbeps_{1}''\sim N(\bsb{0}, \bsb{I}_{d_{z}\times d_{z}})$. Using the standard bound of the normal tail probability, we get
$
P(A^c) \leq 2d_{z}  \varphi(M)/M,  
$ 
where $M=\left(1- \frac{\kappa  d_{nz}}{\mu} \right)\frac{\tau}{\sigma}$.

To bound $P(V^c)$, suppose the spectral decomposition of $\bsbSig_{nz}$ is given by
$\bsbSig_{nz}=\bsb{U} \bsb{D} \bsb{U}^T$ with $\bsb{U}=[\bsb{u}_{1},\cdots,\bsb{u}_{d_{nz}}]^T$, then we can
represent $\bsbSig_{nz}^{-1}$ as $\left[\bsb{u}_i^T \bsb{D}^{-1} \bsb{u}_j \right]_{d_{nz}\times d_{nz}}$,
and $\bsbSig_{nz}^{-1} \bsb{s}$ as $\left[\sum_{j=1}^{d_{nz}} s_j \bsb{u}_i^T \bsb{D}^{-1} \bsb{u}_j
\right]_{d_{nz}\times 1}$. It follows that $\mbox{diag}(\bsbSig_{nz}^{-1})\leq \frac{1}{\mu}$  and
$|\bsbSig_{nz}^{-1} \bsb{s}| \leq \frac{d_{nz}}{\mu}$. Therefore,
$P(V^c)\leq P\left(
\left|\bsbSig_{nz}^{-1} \bsbX_{nz}^{T} \bsbeps \right|\leq L_0\right),   
$ where $L_0=\min|\bsbb_{nz}| - \tau d_{nz} / \mu$.

Because $\bsbSig_{nz}^{-1} \bsbX_{nz}^{T} \bsbeps \sim N(\bsb{0}, \sigma^2 \bsbSig_{nz}^{-1})$, and we have
shown  $\mbox{diag}(\bsbSig_{nz}^{-1}) \leq \frac{1}{\mu}$,  applying Lemma \ref{dawl:lemma2} yields
$$
P(V^c) = P(\max \left|\bsbSig_{nz}^{-1} \bsbX_{nz}^{T} \bsbeps \right| \geq L_0) \leq P(\max |\bsbeps_{2}''|
\geq L_0 \frac{\sqrt \mu}{\sigma}) \leq 2 d_{nz} \Phi([L, \infty)),
$$ where $\bsbeps_{2}'' \sim N(\bsb{0}, \bsb{I}_{d_{nz} \times d_{nz}})$. Hence,
\begin{eqnarray}
1-p_e \leq  P( A^c\cup V^c ) \leq P(A^c) + P(V^c)
 \leq  2 d_{z} \frac{1}{M} \varphi(M) +2 d_{nz} \frac{1}{L} \varphi(L), \label{dawl:(A2.16)}
\end{eqnarray}
where $L
= \frac{\sqrt{\mu }}{\sigma} \left(   \min|\bsbb_{nz}| - \frac{\tau d_{nz}}{ {\mu }}\right)$.

In fact, we can get something slightly stronger than \eqref{dawl:(A2.16)}. Observing that ${\bsbX_{z}'}^T \bsbeps$ is \emph{independent} of $\bsbX_{nz}^T \bsbeps$, we have
\begin{eqnarray}
p_s \geq P(A\cap V) = P(A) \cdot P(V)
\geq \left[ 1 -
2 \Phi (-M)\right]^{d_{z}} \left[ 1 - 2 \Phi (-L)\right]^{d_{nz}}.   \label{dawl:(A2.17)}
\end{eqnarray}
\eqref{dawl:(A2.17)} implies \eqref{dawl:(A2.16)}.

We assumed $\bsbx_i^T\bsbx_i=1$ $i=1,\cdots, d$ in the above derivation. If the $l_2$-norm of each column of $\bsbX$ is no greater than $\sigma_{\max}$, we only need to replace  $\bsbb$, $\hat\bsbb$, $\tau$, by $\bsbb\cdot \sigma_{\max}$, $\hat\bsbb\cdot \sigma_{\max}$, $\tau/\sigma_{\max}$, respectively. The proof of Theorem \ref{selbytisp} is now complete if $\sigma_{\max}=\sqrt n$.\\

For  Theorem \ref{selbyhard}, noticing that (a) $\gsgn(u)=0, \forall |u|>c\tau$ by definition and (b) $p_s \geq P\left(\hat\bsbb_z^{(s)}=0, \mbox{ and } |\hat\bsbb_{nz}^{(s)} |>c\tau^{(s)}\right)$ with $\bsbb^{(s)}=\bsbb\sqrt n$, $\tau^{(s)}=\tau/\sqrt n$,  we can prove it following  the same lines.

\subsection{Proof of Theorem \ref{estbytisp}}
\label{appproofEstByTisp}
Use the same symbols and notations as defined in Appendix \ref{appproof2}.
Let $r_z=\|\hat\bsbb_z\|_2^2$, $r_{nz}=\|\hat\bsbb_{nz} -
\bsbb_{nz}\|_2^2$.
From \eqref{dawl:(A2.9)}, we have
$$
R_{nz} \leq 3 \left( E \| \bsbSig_{nz}^{-1} \bsbX_{nz}^{T} \bsbeps \|_2^2 +
 \tau^2 E\| \bsbSig_{nz}^{-1}\widetilde{{\mbox{sgn}}}(\hat\bsbb_{nz}) \|_2^2
+ \kappa^2 \frac{d_z d_{nz}}{\mu^2} \cdot R_z \right),
$$ due to Cauchy-Schwarz inequality and the fact that $\| \bsbSig_{z, nz} \|_2 = \underset{\|\bsb{\alpha}\|_2=1}{\max} \| \bsbSig_{z, nz} \bsb{\alpha} \|_2 \leq \kappa \sqrt{d_z d_{nz}}$.
Thus
$
R_{nz} \leq 3 \left( \sigma^2 tr(\bsbSig_{nz}^{-1} ) +  \tau^2 \frac{d_{nz}}{\mu^2} +
\kappa^2 \frac{d_z d_{nz}}{\mu^2} \cdot R_z \right)  
$ 
 and \eqref{nzriskbnd} holds.

To get \eqref{zriskbnd}, we need the following result about $\mu_{max}(\bsbS_z^{-1})$, the largest
eigenvalue of $\bsbS_z^{-1}$:
$
\mu_{max}(\bsbS_z^{-1}) \leq \frac{1}{\nu} \left( 1-\kappa^2 \cdot \frac{d_z d_{nz}}{\mu
\nu}\right)^{-1}.
$ 
This is true by noting that $\bsbS_z={\bsbX_z'}^T {\bsbX_z'}$ is semi-positive definite and $\mu_{min}(\bsbS_z) \geq \nu - \kappa^2 d_z d_{nz} / \mu$.

By \eqref{dawl:(A2.8)} and the results in Appendix \ref{appproof2}, we have
\begin{eqnarray*}
r_z &\leq& 3 \mu_{max}^2(\bsbS_z^{-1})  \left[ \|{\bsbX_z'}^{T} \bsbeps\|_2^2 + \tau^2
\|\bsbSig_{z, nz} \bsbSig_{nz}^{-1} \widetilde{{\mbox{sgn}}}(\hat\bsbb_{nz})\|_2^2 + \tau^2
\|\widetilde{{\mbox{sgn}}}(\hat\bsbb_z)\|_2^2 \right] \cdot 1_{A^c} \\
&\leq& 3 \mu_{max}^2(\bsbS_z^{-1})  \left[ \|{\bsbX_z'}^{T} \bsbeps\|_2^2 \cdot 1_{A^c} +
\left(\tau^2 d_z \cdot \kappa^2 \frac{d_{nz}^2}{\mu^2} + \tau^2 d_z\right) \cdot
1_{A^c}\right].
\end{eqnarray*}
Therefore,
\begin{eqnarray*}
R_z &\leq& 3 \mu_{max}^2(\bsbS_z^{-1})  \left[ E(\|{\bsbX_z'}^{T} \bsbeps\|_2^2 \cdot 1_{A^c}) +
\left( 1+ \kappa^2 \frac{d_{nz}^2}{\mu^2} \right) \tau^2 d_z \cdot 2d_z \frac{1}{M}
\varphi(M) \right] \\
&\leq& 3 \mu_{max}^2(\bsbS_z^{-1}) \cdot d_z E\left[( \max|\bsbeps_1'| )^2; \max|\bsbeps_1'| \geq
\tau(1-\kappa\frac{d_{nz}}{\mu}) \right] + \\
& & 3 \mu_{max}^2(\bsbS_z^{-1}) \cdot \frac{1+\kappa^2d_{nz}^2/\mu^2}{(1-\kappa
d_{nz}/\mu)^2}\cdot 2\sigma^2d_z^2 M\varphi(M).
\end{eqnarray*}

Since for a random variable $z$ with probability density $p(\cdot)$ and $a>0$,
$$
E(z^2; z\geq a)=\int_a^{\infty} t^2p(t)dt=\int_a^{\infty} 2s P(z\geq s)ds + a^2 P(z \geq a),
$$ it follows from Lemma \ref{dawl:lemma2} that
$$
E\left[( \max|\bsbeps_1'| )^2; \max|\bsbeps_1'| \geq \tau(1-\kappa\frac{d_{nz}}{\mu}) \right] \leq
E\left[\max|\bsbeps_1''|^2 \cdot \sigma^2; \max|\bsbeps_1''|\geq M \right].
$$ (Recall that $\bsbeps_1''\sim N(\bsb{0}, \bsb{I}_{d_z\times d_z})$.) The density of $\max|\bsbeps_1''|$
is given by $2d_z\varphi(t)(1-2\Phi(-t))^{d_z-1}$. It is easy to get
$$
E\left[\max|\bsbeps_1''|^2 ; \max|\bsbeps_1''|\geq M \right] \leq 2d_z \int_{M}^{\infty}
t^2 \varphi(t)dt \leq 2d_z (M+1/M) \varphi(M).
$$ Hence,
\begin{eqnarray*}
R_z &\leq& 3\frac{1}{\nu^2} \left( 1-\kappa^2 \frac{d_z d_{nz}}{\mu\nu}\right)^{-2} \cdot
\left(
2\sigma^2d_z^2\left(M+\frac{1}{M}\right)\varphi(M)\right.\\
&&\left.+\frac{1+\kappa^2d_{nz}^2/\mu^2}{(1-\kappa
d_{nz}/\mu)^2}\cdot 2\sigma^2d_z^2 M\varphi(M)\right) \\
&=& \frac{\sigma^2}{\nu^2} d_z^2 ( K_1 M + K_2 \frac{1}{M}) \varphi(M).
\end{eqnarray*}

Using a similar scaling argument we obtain Theorem \ref{estbytisp}.

\subsection{Proof of Theorem \ref{orthoestbytisp}}
\label{appproofOrtho}
Let $\hat\bsbb^{H}, \hat\bsbb^{S}$ denote the hard- and soft-thresholding estimates with
threshold value $\tau \sigma$. It is easy to see $\hat\bsbb^{H}$, $\hat\bsbb^{S}$, and
$\hat\bsbb$ all have the same sign and $\hat\bsbb$ is sandwiched by the other two. Therefore,
$
E\| \hat\bsbb - \bsbb \|_2^2 \leq \sum E( \max((\hat\beta_i^{S}-\beta_i)^2,
(\hat\beta_i^{H}-\beta_i)^2)).
$
It is sufficient to study soft- and hard-thresholdings in the univariate case.

Let $y=\mu+\epsilon$ (all are scalars) with $\epsilon \sim N(0, 1)$, and $\rho_S (\tau,\mu)$,
$\rho_H (\tau,\mu)$ be the risks of the soft- and hard-thresholdings with parameter $\tau$.
It is well known~\cite{Donoho, CandOver} that
\begin{eqnarray}
\rho_S(\tau, \mu) \leq \min(\rho_S(\tau, 0)+\mu^2, 1+\tau^2)\leq \min(
\frac{2\varphi(\tau)}{\tau}+\mu^2, 1+\tau^2) \label{dawl:(A4.1)}
\end{eqnarray} for any $\tau >0$. Yet it seems that there is {no} such  explicit nonasymptotic bound, or a
complete proof for the hard-thresholding rule. This short appendix is mainly to give some details
about
this.\\

Our goal is to show the following on the basis of~\cite{Donoho}
\begin{eqnarray}
\rho_H(\tau,\mu) &\leq &1+\tau^2 {\mbox{ for }} \tau > 1  \label{dawl:(A4.2)}
\\
\rho_H(\tau,\mu) &\leq &\rho_H(\tau, 0)+1.2 \mu^2.  \label{dawl:(A4.3)}
\end{eqnarray}
Dohoho \& Johnstone have shown \eqref{dawl:(A4.2)}, and \eqref{dawl:(A4.3)} for $0<\mu<\tau$, but it is technically
difficult  to use the second derivative to prove \eqref{dawl:(A4.3)} for any
$\mu>0$.
Let $g={\partial \rho_H}/{\partial \mu}-2.4\mu$, and $\rho_H(\tau, \mu)$ is known to
be~\cite{Donoho, Gao}
$$
1+(\mu^2-1)(\Phi(\tau-\mu)-\Phi(-\tau-\mu)) + (\tau+\mu) \varphi(\tau+\mu) +
(\tau-\mu)\varphi(\tau-\mu),
$$ where $\varphi, \Phi$ are the standard normal density and distribution functions, respectively.
One may observe that
$
\sup_{\mu\geq 0} g(0, \mu) \leq \sup_{\tau \geq 0} g(\tau, 0) = 0,
$
which is trivial to verify. So it is sufficient to show that for any $(\tau, \mu)>0$, there
exists some $\theta \in [\pi, \frac{3}{2} \pi]$ such that the directional derivative $D_\theta g$ at $(\tau, \mu)$ is greater than $0$, or $\exists \theta_{\tau,\mu} \in [0,\frac{\pi}{2}]$ s.t. $D_\theta g(\tau, \mu) < 0$,
because $g$ is smooth enough.

Consider a uniform direction $\theta=\frac{\pi}{4}$, and let $h=D_\theta g=(\frac{\partial
g}{\partial \tau}+\frac{\partial g}{\partial \mu})/\sqrt 2$. We assume $\mu \geq \tau$ in the
following. Then simple calculations yield
\begin{eqnarray*}
h(\tau,\mu)&=&\sqrt 2 \cdot [(\Phi(\tau + \mu)-\Phi(\mu-\tau)) -\mu \varphi (\mu - \tau) +\\
& & \varphi(\tau+\mu)(\tau^3+3\tau^2\mu+2\tau\mu^2+\mu-2\tau)-1.2 ] \\
&\leq& \sqrt 2 (0.5 + (\tau + \mu)^3 \varphi(\tau + \mu)-1.2) \\
&\leq& (0.5+0.5-1.2) < 0.
\end{eqnarray*}
Therefore,
\begin{equation}
\rho_H(\tau,\mu) \leq \min (\rho_H(\tau, 0)+1.2\mu^2, 1+\tau^2) \leq \min (2
\varphi(\tau)(\tau + \frac{1}{\tau})+1.2\mu^2, 1+\tau^2)   \label{dawl:(A4.4)}
\end{equation} for any $\tau>1$.
Now, combining \eqref{dawl:(A4.1)} and \eqref{dawl:(A4.4)} we can bound the univariate TISP risk
\begin{eqnarray*}
\rho(\tau, \mu) &\leq& \max(\rho_S(\tau, \mu), \rho_H(\tau, \mu)) \leq (1+\tau^2) \min\left(
\frac{2\varphi(\tau)}{\tau}+\frac{1.2}{1+\tau^2}\mu^2, 1 \right) \\
&\leq& (1+\tau^2) \min\left( \frac{2\varphi(\tau)}{\tau}+\mu^2, 1 \right),
\end{eqnarray*}
for any $\tau>1$.  Theorem \ref{orthoestbytisp} thus follows.

Finally it may be worth mentioning that although applying Stein's lemma  is one
possible way (see, for example, Gao~\cite{Gao}), it does not handle  the oracle bound well for an estimator very close to hard
thresholding --- like Zou's oracle bound for the adaptive lasso~\cite{Zou},
because the hard-thresholding function is \emph{not} weakly differentiable.
(Due to an error made in the derivative calculation, Zou's oracle bound for the adaptive lasso defined by $\min \frac{1}{2} \|\bsby - \bsbX \bsbb\|_2^2+\tau \sum w_i |\beta_i| $ with $w_i'\propto|\hat \beta_{ols,i}|^{-\eta}$ should be $(2 \log n + 5 + 4\eta)\cdot \left( \sum \min(\beta_i^2,\sigma^2)+ {\sigma^2}/{(2\sqrt{\pi \log n})}\right),$ with the first factor being $(2 \log n + 5 + 4\eta)$ instead of $(2 \log n + 5 + 4/\eta)$, which \emph{diverges} as $\eta$ goes to infinity. See~\cite{Shethesis} for detail.)

\subsection{Proof of Theorem \ref{hybtispTH}}
\label{appproof3}


In the proof, all inequalities and the absolute value `$||$' are understood
in the componentwise sense.

First we calculate the generalized sign for the hybrid-thresholding \eqref{hybridthfunc}
\begin{eqnarray}
\gsgn (u;\lambda, \eta) = \begin{cases}
            \in [-1,1], \quad\mbox{ if } u = 0\\
            0, \quad\quad\mbox{ if } |u| \in (0, \frac{\lambda}{1+\eta})\\
            \frac{\eta}{\lambda} \cdot u, \quad\mbox{ if } |u|\geq \frac{\lambda}{1+\eta}
            \end{cases}.
\end{eqnarray}
And note that $\tau(\lambda)=\lambda$.
The generalized sign form of the $\Theta$-equation for  Hybrid-TISP estimate $\hat\bsbb$ from \eqref{tispvar} is
\begin{eqnarray}
\bsbSig\hat\bsbb=\bsbX^T \bsby - \lambda \gsgn\left(\hat\bsbb;\frac{\lambda}{k_0^2},\frac{\eta}{k_0^2}\right), \label{kktforhybridtisp}
\end{eqnarray}
where $k_0=\|\bsbX\|_2$.

The proof still follows the lines of the proof for Theorem \ref{selbytisp}. Assume, for the moment, $\bsbX$ has been column-normalized such that the diagonal entries of $\bsbSig=\bsbX^T\bsbX$ are all 1. Clearly, $\hat\bsbb_z=\bsb{0}$, $|\hat\bsbb_{nz}|\geq \frac{\lambda}{k_0^2+\eta}$ is a sufficient condition for the zero consistency of $\hat\bsbb$. From Lemma \ref{dawl:lemma1}, the $\Theta$-equation is equivalent to
\begin{eqnarray*}
\begin{cases}
\bsb{S}_z\hat\bsbb_z=
(\bsbX_z^T-\bsbSig_{z,nz}\bsbSig_{nz}^{-1}\bsbX_{nz}^T)\bsbeps+
\lambda\bsbSig_{z,nz}\bsbSig_{nz}^{-1}\gsgn(\hat\bsbb_{nz};\frac{\lambda}{k_0^2},\frac{\eta}{k_0^2})-\lambda\gsgn(\hat\bsbb_z;\frac{\lambda}{k_0^2},\frac{\eta}{k_0^2})\\
\hat\bsbb_{nz}= \bsbb_{nz}+\bsbSig_{nz}^{-1}(\bsbX_{nz}^T\bsbeps-\lambda\gsgn(\hat\bsbb_{nz};\frac{\lambda}{k_0^2},\frac{\eta}{k_0^2}))-\bsbSig_{nz}^{-1}\bsbSig_{z,nz}\hat\bsbb_z
\end{cases}.
\end{eqnarray*}
Our calculations based on the definition of $\gsgn$ show that
\begin{eqnarray*}
\begin{cases}
\lambda\gsgn(\bsb{0})=\left\{\bsbX_z^T-\bsbSig_{z,nz}\bsbSig_{nz}^{-1}[\bsbI-\eta(\bsbSig_{nz}+\eta\bsbI)^{-1}]\bsbX_{nz}^T\right\}\bsbeps+\eta\bsbSig_{z,nz}(\bsbSig_{nz}+\eta\bsbI)^{-1}\bsbb_{nz}\\
\hat\bsbb_{nz}= (\bsbSig_{nz}+\eta\bsbI)^{-1}\bsbSig_{nz}\bsbb_{nz}+(\bsbSig_{nz}+\eta\bsbI)^{-1}\bsbX_{nz}^T\bsbeps
\end{cases}.
\end{eqnarray*}
Define
\begin{eqnarray*}
A&\triangleq & \left\{\left|\left\{\bsbX_z^T-\bsbSig_{z,nz}\bsbSig_{nz}^{-1}[\bsbI-\eta(\bsbSig_{nz}+\eta\bsbI)^{-1}]\bsbX_{nz}^T\right\}\bsbeps+\eta\bsbSig_{z,nz}(\bsbSig_{nz}+\eta\bsbI)^{-1}\bsbb_{nz}\right|\leq \lambda\right\}\\
V&\triangleq &
\left\{\left|(\bsbSig_{nz}+\eta\bsbI)^{-1}\bsbSig_{nz}\bsbb_{nz}+(\bsbSig_{nz}+\eta\bsbI)^{-1}\bsbX_{nz}^T\bsbeps\right|\geq \frac{\lambda}{k_0^2+\eta}\right\}.
\end{eqnarray*}
Then $p_e\leq P(A^c\cup V^c)\leq P(A^c)+P(V^c)$.

To bound the first probability, noticing that
$
\left|\eta\bsbSig_{z,nz} (\bsbSig_{nz}+\eta\bsbI)^{-1}\bsbb_{nz}\right|\leq  \kappa\sqrt{d_{nz}} \frac{\eta}{\mu+\eta} \|\bsbb_{nz}\|_2,
$ 
we have
$
P(A^c)\leq P\left(\max |\bsbeps_1'|\geq \lambda-\kappa\sqrt{d_{nz}} \frac{\eta}{\mu+\eta} \|\bsbb_{nz}\|_2\right),
$ 
where $\bsbeps_1'=\left\{\bsbX_z^T-\bsbSig_{z,nz}\bsbSig_{nz}^{-1}[\bsbI-\eta(\bsbSig_{nz}+\eta\bsbI)^{-1}]\bsbX_{nz}^T\right\}\bsbeps$.
Since
$$
var(\bsbeps_1')=\sigma^2\left\{\bsbSig_z - \bsbSig_{z,nz}[\bsbI-\eta^2(\bsbSig_{nz}+\eta\bsbI)^{-2}]\bsbSig_{nz}^{-1}\bsbSig_{z,nz}^T\right\},
$$
$\mbox{diag}(var(\bsbeps_1'))\leq \sigma^2\mbox{diag}(\bsbSig_z)\leq \sigma^2 \bsb{1}$. It follows from Lemma \ref{dawl:lemma2} that
$
P(A^c)\leq P\left(\max|\bsbeps_1''|\sigma\geq \lambda - \kappa\sqrt{d_{nz}} \frac{\eta}{\mu+\eta} \|\bsbb_{nz}\|_2\right),
$
where $\bsbeps_1''\sim N(\bsb{0}, \bsbI_{d_{nz}\times d_{nz}})$. Define $M''=\frac{1}{\sigma}\left(\lambda-\kappa\frac{\eta}{\mu+\eta}\sqrt{d_{nz}}\|\bsbb_{nz}\|_2\right)$. We obtain
$
P(A^c)\leq 2 d_z \Phi([M'',+\infty))\leq 2 d_z \varphi(M'')/M''.
$

Next  consider $P(V^c)$. Let $\bsbeps_2'=(\bsbSig_{nz}+\eta\bsbI)^{-1}\bsbX_{nz}^T\bsbeps$. Then
$$
P(V^c)\leq P\left(\max|\bsbeps_2'|\geq \iota-\frac{\lambda}{k_0^2+\eta}\right).
$$
Since $var(\bsbeps_2')=(\bsbSig_{nz}+\eta\bsbI)^{-1} \bsbSig_{nz} (\bsbSig_{nz}+\eta\bsbI)^{-1} \sigma^2$, $\mbox{diag}(var(\bsbeps_{2}'))\leq\frac{\mu\sigma^2}{(\mu+\eta)^2}$. By Lemma \ref{dawl:lemma2} again, we know
$
P(V^c)\leq P\left(\max |\bsbeps_2''|\geq\frac{\mu+\eta}{\sqrt \mu \sigma}\left(\iota-\frac{\lambda}{k_0^2+\eta}\right)\right),
$
where $\bsbeps_2''\sim N(\bsb{0}, \bsbI_{d_{nz}\times d_{nz}})$. Define $L''=\frac{\mu+\eta}{\sqrt \mu \sigma} \left(\iota - \frac{\lambda}{k_0^2+\eta}\right)$. It follows that
$
P(V^c)\leq 2 d_{nz} \Phi([L'',+\infty))\leq 2 d_{nz} \varphi(L'')/L''.
$
\\

We assumed $\bsbx_i^T\bsbx_i=1$ $i=1,\cdots, d$ in the above derivation. If the $l_2$-norm of each column of $\bsbX$ is no greater than $\sigma_{\max}$, it is not difficult to know that we only need to replace the $\bsbb$, $\hat\bsbb$, $\lambda$, $\eta$, by $\bsbb\cdot \sigma_{\max}$, $\hat\bsbb\cdot \sigma_{\max}$, $\lambda/\sigma_{\max}$, $\eta/\sigma_{\max}^2$, respectively. The proof of Theorem \ref{hybtispTH} is now complete if $\sigma_{\max}=\sqrt n$.

\bibliographystyle{acm}
\bibliography{tispbib}

%
%
%
%
%
%
%
%
%

\end{document}